\documentclass[reqno,A4paper]{amsart}

\usepackage{amsmath}
\usepackage{amssymb}
\usepackage{graphicx}
\usepackage{amsthm}
\usepackage{enumerate}
\usepackage{mathrsfs} 
       \usepackage{xcolor}                      
\usepackage{eqlist}
\usepackage{array}
\usepackage{hyperref}
\newtheorem{proposition}{Proposition}

\newtheorem{corollary}[]{Corollary}

\theoremstyle{plain}
\newtheorem{theorem}{Theorem}[section]

\newtheorem{lemma}{Lemma}[section]

\usepackage{url}

\theoremstyle{definition}
\newtheorem{definition}{Definition}
\newtheorem{remark}{\textup{Remark}} 
\newtheorem{example}{\textit{Example}} 

\numberwithin{equation}{section}

\begin{document}

\title[]%
{Some ordering properties of highest and lowest order statistics with exponentiated Gumble type-II distributed components}
\author[Stochastic comparison results for exponentiated Gumbel type-II distribution ]%
{Surojit Biswas}

\newcommand{\acr}{\newline\indent}

\address{\llap{*\,}Department of Mathematics\acr
                   IIT Kharagpur\acr
                   Kharagpur, West Bengal\acr
                   Pin-721302\acr
                   India}
\email{surojit23.iitkgp@gmail.com}





\begin{abstract}

This paper studies stochastic comparisons of the highest and lowest order statistics arising from independent Exponentiated Gumbel Type-II distributed random variables. First, using vector majorization, we establish a comprehensive set of stochastic orderings for the lowest and highest order statistics with respect to the outer and inner shape parameters. In particular, for the lowest order statistic, we derive results on the usual stochastic order and likelihood ratio order, which further imply failure rate, mean residual life, and increasing convex orderings. Second, we show that the lowest order statistics satisfy a proportional failure rate structure. In addition, we also establish star and dispersive orderings.
For the highest order statistic, we obtain results in the reversed failure rate order, and we show that, in general, a proportional failure rate structure does not hold. Moreover, unlike several classical exponentiated models, the highest order statistics are not necessarily comparable in the usual stochastic order under heterogeneity of the inner shape parameters. 
Finally, we extend the analysis using multivariate chain majorization to study the effect of simultaneous variation in scale and outer shape parameters. Under this framework, we establish usual stochastic ordering results for both the lowest and highest order statistics. It provides more general reliability comparisons under joint parameter heterogeneity.
\end{abstract}

\maketitle

\section{Introduction}

In nature, parallel and series systems are very common complex systems to perform reliability optimization. A significant amount of research has been conducted on reliability optimization and stochastic ordering due to their importance in various types of complex systems (complex networks). One may refer to \cite{kuo2000annotated} for this kind of optimization problem. Parallel and series systems are $1$-out-of-$n$ and $n$-out-of-$n$ systems, respectively. A $k$-out-of-$n$ system is a system that functions if and only if at least $k$ out of its $n$ components function.
These systems fall under a particular class of order statistics. The study of stochastic comparisons of $k$-among-$n$ systems is very similar to that of stochastic comparisons of order statistics. \\

Stochastic comparisons of order statistics play an essential role in many
fields, such as statistics, applied probability, reliability theory, radio communication, auction
theory, operation management, etc. For example, the application of stochastic ordering in wireless communication and operation management can be found in the article by \cite{tepedelenlioglu2011applications} and \cite{feng2022applications}, respectively.  Within these fields, particularly in reliability and engineering, these comparisons help answer crucial questions for a quality control manager or design engineer. For example, how does one compare two systems with different combinations of component characteristics? Is a system composed of components with highly variable strength parameters stochastically weaker than a system with more uniform components, even if the average strength is the same? This question is similar to understanding that a chain depends on its weakest link. This paper attempts to answer exactly these kinds of questions by offering rigorous methods to compare the reliability of series (lowest order statistic) and parallel (highest order statistic) systems using potent vector and multivariate majorization techniques. These methods provide useful tools for reliability optimization and system design decision-making.\\

  In this article, we discuss various results on stochastic ordering for the parallel and series systems where the components of the system come from the independent exponentiated Gumbel type-II distributions. A random variable $ X $ is said to have an exponentiated Gumbel type-II distribution (in short, we use $EG_{2}$ this paper) due to \cite{okorie2016exponentiated}  having the following cumulative distribution function (cdf) 

\begin{equation}
 F(x)= 1-\left( 1-e^{-\theta x^{-\phi}}\right)^{\alpha} ,\hspace{.8cm} x> 0; \hspace{.25cm}\theta, \phi, \alpha >0.
 \label{eq1}
 \end{equation}

   Here $ \theta $ is the scale parameter, and $ \phi, \alpha $ are shape parameters. We call  $ \phi $ as the inner shape parameter and $ \alpha$ the outer shape parameter. We use the notation $ X\sim EG_{2}(\theta,\phi,\alpha) $ if $ X $ has the cumulative distribution function in equation (\ref{eq1}).
  This $ EG_{2} $ distribution is a new lifetime distribution.  By simply assigning the values of $ \alpha = 1, \theta = 1 $, and $ \alpha, \theta = 1 ,$ one can have the following distributions:  Gumbel type-2 distribution,
 Exponentiated Fr\'echet (EF) distribution, and Fr\'echet distribution from this $ EG_{2}$ distribution, respectively, and for $ y=x^{-\phi} $, it becomes
 Exponentiated Exponential (EE) distribution.


The Exponentiated Gumbel Type-II ($EG_{2}$) distribution provides a flexible framework for modeling the lifetimes of components exhibiting heavy-tailed and extreme-value behavior, which frequently arises in reliability engineering applications such as material strength, fatigue failure, and environmental stress modeling. 
From a reliability perspective, the parameters of the $EG_{2}$ distribution admit a natural physical interpretation. The inner shape parameter $\varphi$ governs the intrinsic aging characteristics and component-level heterogeneity arising from variations in material properties or manufacturing imperfections. In contrast, the outer shape parameter $\alpha$ acts as a global transformation parameter that modifies the baseline failure distribution. This mechanism allows the model to capture system-wide amplification or attenuation of failure risk, which may arise due to common environmental effects or shared stress-loading conditions.
Such a separation between local (component-specific) and global (system-level) effects is particularly relevant in complex engineering systems. Moreover, since the lifetime of a series system is determined by the lowest order statistic and that of a parallel system by the highest order statistic, the $EG_{2}$ model provides a mathematically tractable and physically grounded setting for studying system reliability through stochastic comparisons of order statistics under parameter heterogeneity. 
In addition, Okorie et al. \cite{okorie2016exponentiated} demonstrated, based on log-likelihood and information criteria, that the $EG_{2}$ distribution provides a superior fit to real-life data compared to several classical lifetime models, including the Gumbel Type-II, Fr\'echet, exponentiated Fr\'echet, Weibull, and log-normal distributions. This empirical flexibility, combined with its clear reliability interpretation, makes the $EG_{2}$ distribution a suitable and practically relevant choice for modeling component lifetimes and analyzing system reliability. Now, let $ X_{1}, X_{2},\cdots,X_{n} $ be a set independently distributed random variables and \\$$ X_{n:n}=\max\left\lbrace X_{1},\cdots,X_{n}\right\rbrace,  ~~~~X_{1:n}=\min\left\lbrace X_{1}, X_{2},\cdots,X_{n}\right\rbrace.$$ 
 $ X_{1:n} $ is known as the $1$st (lowest) order statistic, which represents a series system, and $ X_{n:n} $ is known as the $ n $th (highest) order statistic, which represents a parallel system. For more detailed information about order statistics, see \cite{shaked2007stochastic}. Various discussions on order statistics in terms of stochastic comparisons are already available where the component variables follow  Generalized Exponential distribution (see, \cite{balakrishnan2014stochastic}), Exponential Weibull distribution (see, \cite{fang2016ordering}), Exponentiated Scale model distribution (see, \cite{bashkar2017stochastic}), Fr\'echetdistribution  (see \cite{gupta2015stochastic}),  Gumbel distribution  (see, \cite{biswas2022stochastic}), Marshall and Olkin’s Family of Distributions (see \cite{das2021some}), Exponentiated Location-Scale models (see \cite{das2021ordering})\textit{etc.} The comparisons include usual stochastic order, failure rate order, reversed failure rate order, likelihood ratio order, dispersive order, \textit{etc}. For further details on stochastic
 comparisons, one may refer to  \cite{misra2012new},  \cite{torrado2015stochastic},  \cite{fang2012new}, \cite{toomaj2017some}, and  \cite{nadarajah2006exponentiated}.
 Some results from  \cite{balakrishnan2014stochastic} are helpful in developing this article.  More recently, dispersive, star, and related variability orders have been investigated for extreme order statistics and complex reliability structures. For example, dispersive and star orders for the lowest and highest order statistics have been established for dependent unit gamma Gompertz random variables under dependence structures such as Archimedean copulas (see \cite{shekari2025dispersion}), while dispersive ordering of fail‑safe and k‑out‑of‑n systems with heterogeneous exponential components has also been studied (see \cite{zhao2011dispersive}).


 Let $X_{1},\cdots,X_{n}$ be independent random variables such that 
$X_{i}\sim EG_{2}(\theta_{i}, \phi_{i}, \alpha_{i}), \\ i=1,2,\ldots,n$. 
Similarly, let $X^{*}_{1},\cdots,X^{*}_{n}$ be another set of independent random variables with 
$X^{*}_{i}\sim EG_{2}(\theta^{*}_{i}, \phi^{*}_{i}, \alpha^{*}_{i}), \hspace{.2cm} i=1,2,\ldots,n$.
First, we consider the case where the scale and inner shape parameters are common, i.e., 
$\theta_{1}=\cdots=\theta_{n}=\theta^{*}_{1}=\cdots=\theta^{*}_{n}=\theta$ and 
$\phi_{1}=\cdots=\phi_{n}=\phi^{*}_{1}=\cdots=\phi^{*}_{n}=\phi$, and the outer shape parameter vectors satisfy 
$(\alpha_{1},\cdots,\alpha_{n})\succeq^{m}(\alpha^{*}_{1},\cdots,\alpha^{*}_{n})$. 
Under this setup, we establish the usual stochastic order for the lowest order statistic, $X_{1:n}$. We further derive a necessary and sufficient condition for the likelihood ratio order, which consequently implies ordering in terms of the failure rate, mean residual life, and increasing convex orders. Additionally, we establish star and dispersive orderings, providing deeper insights into the variability and aging properties of the corresponding series system. 
For the highest order statistic, $X_{n:n}$, we derive results for the reversed failure rate order. Under the setting of a homogeneous outer shape parameter, we further obtain star and dispersive orderings, which are particularly relevant for parallel system reliability analysis.
Next, we consider the case where the scale and outer shape parameters are common, i.e., $\theta_{i} = \theta_{i}^{*} = \theta$ and $\alpha_{i} = \alpha_{i}^{*} = \alpha$ for all $i$, while the inner shape parameter vectors satisfy $(\phi_{1}, \dots, \phi_{n}) \succeq^{m} (\phi^{*}_{1}, \dots, \phi^{*}_{n})$. In this context, we investigate the effect of heterogeneity in the inner shape parameters and establish usual stochastic ordering for the lowest order statistic, which subsequently implies the increasing convex ordering. Finally, we demonstrate that the highest order statistics are not necessarily comparable under the usual stochastic order when the inner shape parameters are heterogeneous, which highlights a unique non-comparability phenomenon in the $EG_2$ model.

Finally, we extend the analysis to a more general setting using multivariate chain majorization. 
Specifically, when $\phi_{1}=\cdots=\phi_{n}=\phi^{*}_{1}=\cdots=\phi^{*}_{n}$, we compare two systems whose parameter matrices 
$(\theta_{i},\alpha_{i})$ and $(\theta^{*}_{i},\alpha^{*}_{i})$, $i=1,2,\ldots,n$, are related through multivariate majorization. 
Under this framework, we establish usual stochastic ordering results for both the lowest and highest order statistics.

Although majorization techniques and Schur-convexity have been widely used in stochastic ordering of exponentiated families such as exponentiated exponential, exponentiated Weibull, and related models, the present study reveals that the three-parameter $EG_2$ distribution exhibits fundamentally different structural behavior. In particular, unlike simpler models where stochastic comparisons of order statistics often reduce to conditions on aggregate parameters, the $EG_2$ model does not generally admit such simplifications.
We show that, for the highest order statistic, usual stochastic ordering may fail under majorization of the inner shape parameters, leading to a non-comparability phenomenon that does not arise in many classical exponentiated models.
In addition, by applying multivariate chain majorization, we obtain more general stochastic comparisons that capture simultaneous variation in multiple parameters. These results provide new insights into the role of parameter heterogeneity in reliability systems and demonstrate that the $EG_2$ model requires a more refined analytical treatment compared to previously studied families.

 The organization of this article is as follows. 
Section \ref{sec2} gives some basic useful definitions, lemmas, and theorems, which we have used throughout this paper. Section \ref{sec3} has two subsections; in Subsection \ref{subsec1},  we work with the concept of vector majorization technique for some ordering results between lowest and highest order statistics, and in Subsection \ref{subsec2},  we deal with the concept of multivariate majorization and achieve the usual stochastic ordering only. Finally, in Section \ref{con}, we put the conclusion of our discussion.

\section{Preliminaries}
  \label{sec2}
 In this section, we present some well-known definitions together with some results that are most pertinent to the developments of Section \ref{sec3}.  We use the notation  $\mathbb{R}$=$(-\infty,+\infty),$  $ \mathbb{R}_{+}$=$[0,+\infty),$ and $log$ for usual logarithm base $ e,$ throughout this paper.\par
 
 Let $ X, Y $ be two non-negative continuous univariate random variables having the following characteristics
 \begin{itemize}
 \item Cumulative distribution functions: $F(x)$,\hspace{.2cm} $G(x)$.
 \item  Reliability (Survival) functions: $\overline F(x)=1-F(x)$,  \hspace{.2cm} $\overline G(x)=1-G(x)$.

 \item  Probability density functions: $f(x)$,  \hspace{.2cm} $g(x)$.\
 \item Failure (Hazard) rate functions: $r(x)= \dfrac{f(x)}{\overline{F}(x)}, \hspace{.2cm}     s(x)=\dfrac{g(x)}{\overline{G}(x)}$.
 \item Reversed failure rate functions: $\tilde{r}(x)= \dfrac{f(x)}{F(x)}, \hspace{.2cm} \tilde{s}(x)= \dfrac{g(x)}{G(x)}$.
 \end{itemize}

  \begin{definition}
 Let $X, Y$ be two random variables.
 \begin{enumerate}
 \item[1.] $X \leq_{st} Y$ (\textit{stochastic order}) iff $\overline{F}(x) \leq \overline{G}(x)$ $\forall$ $x \in \mathbb{R}$.
 \item[2.] $X \leq_{fr} Y$ (\textit{failure rate order}) iff $s(x) \leq r(x)$ $\forall$ $x \in \mathbb{R}$.
 \item[3.] $X \leq_{rf} Y$ (\textit{reversed failure rate order}) iff $\tilde{r}(x) \leq \tilde{s}(x)$ $\forall$ $x \in \mathbb{R}$.
 \item[4.] $X \leq_{lr} Y$ (\textit{likelihood ratio order}) iff $g(x)/f(x)$ is non-decreasing in $x$.
 \item[5.] $X \leq_{mrl} Y$ (\textit{mean residual life order}) iff $m_X(x) \leq m_Y(x)$ $\forall$ $x \geq 0$, where $m_X(x) = E[X-x|X>x]$.
 \item[6.] $X \leq_{disp} Y$ (\textit{dispersive order}) iff $F^{-1}(\beta) - F^{-1}(\alpha) \leq G^{-1}(\beta) - G^{-1}(\alpha)$ for all $0 < \alpha \leq \beta < 1$.
 \item[7.] $X \leq_{*} Y$ (\textit{star order}) iff $G^{-1}F(x)/x$ is non-decreasing in $x > 0$.
 \item[9.] $X \leq_{icx} Y$ (\textit{increasing convex order}) iff $E[\phi(X)] \leq E[\phi(Y)]$ for all increasing convex functions $\phi$.
 \item[10.] $X \leq_{c} Y$ (\textit{relative failure rate / relative aging order}) iff $r(x)/s(x)$ is non-decreasing in $x$.
 \end{enumerate}
 \end{definition}

The well-known relations between the above definitions are given by:
\begin{equation*}
    X \leq_{lr} Y \implies X \leq_{fr} Y  \implies X \leq_{st} Y,
\end{equation*}
and
\begin{equation*}
    X \leq_{lr} Y \implies X \leq_{rf} Y \implies X \leq_{st} Y.
\end{equation*}

 \begin{definition}

 Let $ \textbf{y}=(y_{1},y_{2},...,y_{n})$ with the order components, $ y_{(1)}\leq....\leq y_{(n)} $ and  $ \textbf{x}=(x_{1},x_{2},...,x_{n})$ with the order components $x_{(1)}\leq...\leq x_{(n)} $, be two real vectors from $ \mathbb{R}^{n} $. Then we say $ \textbf{y}$ majorizes$\textbf{ x} $ if $$ \sum_{i=1}^{k}y_{(i)} \leq \sum_{i=1}^{k}x_{(i)}$$ \hspace{.3cm}$ k=1,2,...,n-1,$ and $\sum_{i=1}^{n}y_{i}= \sum_{i=1}^{n}x_{i}.$ It is denoted by
 $ \textbf{y}\succeq^{m} \textbf{x}$. 
 \end{definition}

 \begin{definition}
 
 Let $ \textbf{y}=(y_{1},\cdots,y_{n})$ and $ \textbf{x}=(x_{1},\cdots,x_{n})$ be two vectors from $ \mathbb{R}^{n} $. A real valued function $ \sigma(\textbf{y}):\mathbb{R}^{n}\rightarrow \mathbb{R} $ is said to be \textit{Schur-convex} and \textit{Schur-concave}  if $ \sigma(\textbf{y})\geq \sigma(\textbf{x})$ and $\sigma(\textbf{y})\leq \sigma(\textbf{x}), \hspace{.2cm} \forall \textbf{y}\succeq^{m} \textbf{x}$, respectively.
 \end{definition}\par
 The theorem stated below is very useful for our results.
 
 \begin{theorem} (\cite{marshall1979inequalities}):\\
 Let $ I\subset\mathbb{R} $ be an
 open interval and let $ \sigma: I^{n} \rightarrow \mathbb{R}$ be continuously differentiable function. The if and only if (iff) conditions for $ \sigma $ to be Schur-convex(Schur-concave) on $ I^{n}$ are
 $ \sigma $ is symmetric on $ I^{n}$ and, for all $ i\ne j $
 \begin{center}
 $ (z_{i}-z_{j})\left( \dfrac{\partial \sigma}{\partial z_{i}}(\textbf{z})-\dfrac{\partial \sigma}{\partial z_{j}}(\textbf{z})\right)\geq 0 \hspace{.1cm} (\leq 0) $
 \end{center}for all $ \textbf{z}\in I^{n} $. Where, $  \dfrac{\partial \sigma}{\partial z_{i}}$ is partial derivative of $ \sigma $ with respect to the $ i^{th} $ component of $ \textbf{z} $.
 
 \label{th1}
 \end{theorem}\par
 
 A square matrix $ \Pi  $ is said to be a permutation matrix if each row and column has only one $1,$ and all other entries are zero. We can always find $n!$ such matrices
 by interchanging rows (or columns) of the identity matrix of order $ n $. Let $ P=(p_{ij}) $ be a matrix of order $n$, $P$ is said to be doubly stochastic if $ p_{ij}\geq 0$ and $\sum_{i=1}^{n} p_{ij}=1=\sum_{j=1}^{n} p_{ij}$ for $ i,j=1,2,....,n $.
 The $ T $-transform matrix has the following form $ T_{w}=wI_{n}+(1-w)\Pi_{n} $, where $ w\in [0,1] $, $ I_{n} $ is the $ n \times n $ identity matrix and $ \Pi_{n} $ is a permutation matrix of order $ n $ that just interchanges two coordinates. One important fact about $ T $-transformation matrices is that the product of a finite number of $ T $-transformation matrices with the same structure is also a $ T $-transformation matrix, and the resulting matrix has the same structure as the elements. But it may not hold for $ T $-transformation matrices with different structures.\\
 
 Now, the next definition deals with various types of multivariate majorization  (see, \cite{marshall1979inequalities}).
 \begin{definition}

 Let $ A=\{a_{ij}\}, B=\{  b_{ij} \}$ be two matrices of size $ m\times n $ such that $ a^{R}_{1},....a^{R}_{m} $ and $ b^{R}_{1},....,b^{R}_{m} $ are the rows of $ A $ and $ B $ respectively, then:
 \begin{enumerate}
 \item[1.] $ B $ is said to be chain majorized by $ A $, denoted by $ A\gg B $ if there exists a finite set of $ T $-transformation matrices $ T_{w_{1}},....,T_{w_{k}} $ of size $ n\times n $ such that $ B=AT_{w_{1}}\cdots T_{w_{k}}$;
 \item[2.] $ B $ is said to be majorized by $ A $, denoted by $ A> B $ if there exists a doubly stochastic matrix $ P $ of size $ n \times n $ such that $ B=AP $;
 \item[3.]  $ B $ is said to be row majorized by $ A $, denoted by $ A>^{row} B $ if $ a^{R}_{i}\succeq^{m} b^{R}_{i} $ for $ i=1,2,....,m$.
 \end{enumerate}
 \label{df4}
 \end{definition}\par
 
 A widespread implementation of these above forms of multivariate majorization is that
 \begin{center}
 $ A\gg B  \implies  A> B \implies  A>^{row} B$.
 \end{center}
 Interested readers may refer to the \cite{marshall1979inequalities} for more details.\par
 
To prove our main results in the next section, we shall use the following theorems. Let us consider three sets $ \mathcal{S}_{n} $, $ \mathcal{T}_{n} $ and $ \mathcal{D}_{n} $ defined as follows
 
 \[  \mathcal{S}_{n}=\left \{
 \begin{bmatrix}
 x_{1},\dots,x_{n} \\
 
 y_{1},\dots ,y_{n}
 \end{bmatrix}:(x_{i}-x_{j})(y_{i}-y_{j})\leq 0,and \hspace{.3cm}x_{i}>0,y_{i}>0,\hspace{.3cm} i,j=1,...,n
 \right \},
 \]
 
 \[\mathcal{T}_{n}=
  \left \{
 \begin{bmatrix}
 x_{1},\dots,x_{n} \\
 
 y_{1},\dots ,y_{n}
 \end{bmatrix}:(x_{i}-x_{j})(y_{i}-y_{j})\leq 0,and \hspace{.3cm}x_{i}\geq 1,y_{i}>0,\hspace{.3cm} i,j=1,...,n
 \right \},
 \]and
 $\mathcal{D}_{n}=\left \{ (x_{1},\dots,x_{n})\in \mathbb{R}^{n}: x_{1}\leq x_{2}\leq \dots \leq x_{n} \right \}.$
 
 \begin{theorem}
  A differentiable function $\psi:\mathbb{R}_{+}^{4}\rightarrow \mathbb{R}_{+} $ satisfies 
 \begin{center}
 $ \psi(A)\leq (\geq) \psi(B) $ for all $ A, B $ such that $ A\in \mathcal{S}_{2}(\mathcal{T}_{2}) $, and $ A\gg B $ if and only if
 \begin{enumerate}
 \item[1.] $ \psi(A)=\psi(\Pi A) $ for all permutation matrices $ \Pi $, and for all  $ A\in \mathcal{S}_{2}(\mathcal{T}_{2}) $;
 \item[2.] $ \sum_{i=1}^{2}(a_{ik}-a_{ij})(\psi_{ik}(A)-\psi_{ij}(A))\leq(\geq) 0 $ for all $ j,k=1,2$ and for all  $ A\in \mathcal{S}_{2}(\mathcal{T}_{2})$, where $ \psi_{ij}(A)=\dfrac{\partial \psi(A)}{\partial a_{ij}} .$
 
 \end{enumerate}
 \end{center}
 \begin{proof}
  Proof of this theorem can be found in \cite{marshall1979inequalities}.
 \end{proof}
 
 \label{th2}
 \end{theorem}
 
 \begin{lemma}
 Let the function $ \eta:(0,\infty)\times(0,1)\rightarrow
 (-\infty,0)$ be defined as
 \begin{center}
 $$ \eta(\alpha,u)=\dfrac{u^{\alpha} \log u}{1-u^{\alpha}}. $$
 \end{center}
 Then, 
 \begin{enumerate}
 \item[1.] $ \eta(\alpha,u) $ is increasing with respect to $ \alpha $ for each $ 0<u<1 $;
 \item [2.] $ \eta(\alpha,u) $ is decreasing with respect to $ u $ for each $ \alpha>0 $.
 \end{enumerate}
 
 \begin{proof}
   For the proof of this lemma, one can see in  \cite{balakrishnan2014stochastic}.
 \end{proof}
 
 \label{lemma3}
 \end{lemma}

 \begin{lemma}
 Let the function $ \gamma:(0,\infty)\times(0,1)\rightarrow
 (0,\infty)$ be defined as 
 \begin{center}
 $$ \gamma(\alpha,u)=\dfrac{\alpha (1-u)u^{\alpha-1}}{1-u^{\alpha}}. $$.
 \end{center}
 Then,
 \begin{enumerate}
 \item[1.] $ \gamma(\alpha,u) $ is decreasing with respect to $ \alpha $ for each $ 0<u<1 $;
 \item [2.] $ \gamma(\alpha,u) $ is decreasing with respect to $ u $ for each $ 0<\alpha \leq 1 ;$
 \item [3.] $ \gamma(\alpha,u) $ is increasing with respect to $ u $ for each $ \alpha \geq 1 $.
 \end{enumerate}
 \begin{proof}
 Proof of this lemma can be found in \cite{balakrishnan2014stochastic}.
 \end{proof}
 
 \label{lemma4}
 \end{lemma}

 \begin{lemma}
Let the function $ \phi(\alpha,u):(0,\infty)\times(0,1)\rightarrow (0,\infty)$ be defined as 
 \begin{center}
 $$ \phi(\alpha,u)=\dfrac{\alpha   u^{\alpha-1}}{1-u^{\alpha}} .$$
 \end{center}
 Then $ \phi(\alpha,u)$ is convex in $ \alpha $ for any  $0<u<1. $
 \begin{proof}
  For the proof in detail, see \cite{balakrishnan2014stochastic}.
 \end{proof}\
 \label{lemma5}
 \end{lemma}

 The following lemma has been taken from book in \cite{marshall1979inequalities}, which will be used to prove one of our results in the next section.
 
 \begin{lemma}
Let $ \mathcal{I}\subset \mathbb{R} $ be an open interval such that if a function $ \lambda(x):\mathcal{I}\rightarrow\mathbb{R}
$ is convex then the function $g(\overline{x}) $ defined as
   $$ g(\overline{x})=\sum_{i=1}^{n}\lambda(x_{i}) $$
is a Schur-convex function in $ \mathcal{I}^{n}, $ where $ \overline{x}=(x_{1},...,x_{n}). $ Consequently for $ \overline{x}=(x_{1},...,x_{n}), $ $ \overline{y}=(y_{1},...,y_{n})\in \mathcal{I}^{n} $ if $ (x_{1},...,x_{n})\succeq^{m}(y_{1},...,y_{n})$ then $ g(\overline{x})\geq g(\overline{y}). $

\label{lemma6}
\end{lemma}

\section{Main results and discussion}
\label{sec3}
\subsection{Results based on vector majorization.}
\label{subsec1}
 This section presents several results by comparing the highest and lowest-order statistics of the exponentiated Gumbel type-II distribution. We put a sufficient condition on the likelihood ratio of the lowest order statistics based on the vector majorization technique for the outer shape parameter $ \alpha,$ and finally, we compare the lowest order statistics by using the usual stochastic order of the inner shape parameter $ \phi. $\\

The lowest order statistic represents the lifetime of a series system, since such a system fails as soon as its first component fails. Therefore, comparing the lowest order statistics is of direct importance in reliability analysis. The following result studies how heterogeneity in the \textit{outer shape parameter} $ \alpha$ (here, $\theta$ is the scale parameter, and $\phi$ is the inner shape
parameter) vector affects the reliability of two series systems in the sense of the usual stochastic order. Such comparisons are useful in identifying conditions under which one system is more reliable than another when the component lifetimes follow the $EG_2$ model.

 \begin{theorem}
 Let $ X_{1},\cdots.,X_{n} ,$ and $ X^{*}_{1},\cdots,X^{*}_{n} $  be two sets of non-negative independent random variables such that $ X_{i}\sim EG_{2}(\theta_{i}, \phi, \alpha_{i}),$ and  $ X^{*}_{i}\sim EG_{2}(\theta_{i}, \phi, \alpha_{i}^{*}),$\\
$ i=1,2,....,n$. Then, if  $ (\alpha_{1},\cdots\alpha_{n})\succeq^{m}(\alpha^{*}_{1},\cdots\alpha^{*}_{n}) $ and $\overline{\theta}=(\theta_{1},\theta_{2},\dots,\theta_{n})\in \mathcal{D}_{n}$, we have\\ $ X_{1:n}\leq_{st} X_{1:n}^{*} $.
  
\begin{proof}
The the reliability function of $ X_{1:n} $
$$\overline{F}_{X_{1:n}}(x)=\prod_{i=1}^{n}\left( 1-e^{-\theta_{i} x^{-\phi}}\right)^{\alpha_{i}}\hspace{.4cm} x,\theta_{i},\phi,\alpha_{i}>0, i=1,...,n. $$
Now, the partial derivative of $\overline{F}_{X_{1:n}}(x) $ with respect to $ \alpha_{i}, i=1,...,n, $ is 
\begin{center}
$ \dfrac{\partial \overline{F}_{X_{1:n}}(x)}{\partial \alpha_{i}}=\overline{F}_{X_{1:n}}(x) \log  \left( 1-e^{-\theta_{i} x^{-\phi}}\right).$

\end{center}
Consider, $ \sigma_{X_{1:n}}(\alpha_{i},\alpha_{j})=(\alpha_{i}-\alpha_{j})\left( \dfrac{\partial \overline{F}_{X_{1:n}}(x)}{\partial \alpha_{i}}-\dfrac{\partial \overline{F}_{X_{1:n}}(x)}{\partial \alpha_{j}}\right) $. $$\sigma_{X_{1:n}}(\alpha_{i},\alpha_{j})=(\alpha_{i}-\alpha_{j})\overline{F}_{X_{1:n}}(x) \left[ \log  \left( 1-e^{-\theta_{i} x^{-\phi}}\right)- \log  \left( 1-e^{-\theta_{j} x^{-\phi}}\right)\right].$$
\\
Without loss of generality let us assume that $\alpha_{i}\geq\alpha_{j}$ and since $ \overline\theta\in \mathcal{D}, $ this implies that $\theta_{j}\geq \theta_{i}, \forall j\geq i .$
 Now, it is clear that $ \log  \left( 1-e^{-\theta x^{-\phi}}\right) $ is increasing with respect to $ \theta .$ So, we have  $$ \log  \left( 1-e^{-\theta_{j} x^{-\phi}}\right)\geq   \log  \left( 1-e^{-\theta_{i} x^{-\phi}}\right), \forall \theta_{j}\geq \theta_{i}. $$ 
 
 Therefore, considering the above-mentioned fact together with the assumption $\alpha_{i}\geq\alpha_{j}$, we conclude that
 $ \sigma_{X_{1:n}}(\alpha_{i},\alpha_{j})\leq 0, $ which implies that $\overline{F}_{X_{1:n}}(x)  $ is \textit{Schur-concave} function with respect to $ (\alpha_{1},\cdotp\cdotp\cdotp,\alpha_{n}).$ So, $\overline{F}_{X^{*}_{1:n}}(x)\geq\overline{F}_{X_{1:n}}(x)  $ i.e.,  $ X_{1:n}\leq_{st} X_{1:n}^{*},$ as required.  
\end{proof} 
\label{th7}
\end{theorem}

\begin{corollary}
       Let $ X_{1},\cdots.,X_{n} ,$ and $ X^{*}_{1},\cdots,X^{*}_{n} $  be two sets of non-negative independent random variables such that $ X_{i}\sim EG_{2}(\theta_{i}, \phi, \alpha_{i}),$ and  $ X^{*}_{i}\sim EG_{2}(\theta_{i}, \phi, \alpha_{i}^{*}),$ 
$ i=1,2,....,n$. Then, if  $ (\alpha_{1},\cdots\alpha_{n})\succeq^{m}(\alpha^{*}_{1},\cdots\alpha^{*}_{n}) $ and $\overline{\theta}=(\theta_{1},\theta_{2},\dots,\theta_{n})\in \mathcal{D}_{n}$, we have $ X_{1:n}\leq_{icx} X_{1:n}^{*} $.
\end{corollary}

\begin{proof}
        The proof follows from the chain of implications existing between usual stochastic orders and the increasing convex orders. Let $U$ and $V$ be two random variables. It is well-established in reliability theory (see 
        Theorem-4.A.55 pp-205, \cite{shaked2007stochastic}) that  $U \leq_{st} V \implies U \leq_{ss} V \implies U \leq_{icx} V$. Hence, the corollary follows.
\end{proof}

The highest order statistic corresponds to the lifetime of a parallel system, because the system continues to operate until all components fail. In reliability theory, the reversed failure rate order is particularly relevant for parallel systems, since it compares their behaviour from the left tail and provides useful information on early-life performance and inactivity risk. The following theorem, therefore, gives a meaningful comparison between two parallel systems through majorization of the \textit{outer-shape} parameter vector.

\begin{theorem}
 Let $ X_{1},\cdots.,X_{n} ,$ and $ X^{*}_{1},\cdots,X^{*}_{n} $  be two sets of non-negative independent random variables such that $ X_{i}\sim EG_{2}(\theta, \phi, \alpha_{i}),$ and  $ X^{*}_{i}\sim EG_{2}(\theta, \phi, \alpha_{i}^{*})$,\\
$ i=1,2,....,n$. Then, if  $ (\alpha_{1},\cdots\alpha_{n})\succeq^{m}(\alpha^{*}_{1},\cdots\alpha^{*}_{n}) $, we have
$ X_{n:n}\geq_{rf} X_{n:n}^{*} $.

\begin{proof}
We know that, in a parallel system, the reversed failure rate of the system lifetime is proportional to the reversed failure rate of the lifetime of each component of the system. Therefore for $ x> 0, $ the reversed failure rate is 
$$\tilde{r}_{X_{n:n}}(x)=\theta \phi x^{-\phi-1} e^{-\theta x^{-\phi}} \sum_{i=1}^{n}\dfrac{\alpha_{i} \left(  1-e^{-\theta x^{-\phi}}\right)^{\alpha_{i}-1}}{1-\left( 1-e^{-\theta x^{-\phi}}\right)^{\alpha_{i}}}, i=1,...,n.$$
Let $ t=1-e^{-\theta x^{-\phi}} $ and it is clear that $ 0<t<1 ,$ for $ x,\theta,\phi,\alpha>0 $. Then
$$ \tilde{r}_{X_{n:n}}(x)=\theta \phi x^{-\phi-1} e^{-\theta x^{-\phi}}\sum_{i=1}^{n}\varphi(\alpha_{i},1-e^{-\theta x^{-\phi}}), $$ where $\varphi(\alpha,t)=\dfrac{\alpha   t^{\alpha-1}}{1-t^{\alpha}}. $ Now, using Lemma \ref{lemma5}, and Lemma \ref{lemma6} we can conclude that $ \sum_{i=1}^{n} \varphi(\alpha_{i},1-e^{-\theta x^{-\phi}}) $ is Schur-convex in $ (\alpha_{1},\cdots,\alpha_{n} )$. Hence, the theorem follows.
\end{proof}

\label{th8}
\end{theorem}








Beyond usual stochastic ordering, the likelihood ratio order is one of the strongest and most informative stochastic orders used in reliability theory. Once established, it automatically implies several other useful orderings, such as failure rate and usual stochastic orderings. For a series system, such a result gives a sharper comparison of system lifetimes and can be used to assess which configuration is uniformly more reliable. Next, we present a \textit{necessary and sufficient condition} for the likelihood ratio ordering of the lowest order statistics.
\begin{theorem}
Let $ X_{1},\cdots.,X_{n} ,$ and $ X^{*}_{1},\cdots,X^{*}_{n} $  be two sets of non-negative independent random variables such that $ X_{i}\sim EG_{2}(\theta, \phi, \alpha_{i}),$ and  $ X^{*}_{i}\sim EG_{2}(\theta, \phi, \alpha_{i}^{*})$,
$ i=1,2,....,n$. Then \\ $ \sum_{i=1}^{n}\alpha_{i}\leq \sum_{i=1}^{n}\alpha^{*}_{i}$  $ \iff X^{*}_{1:n}\leq_{lr} X_{1:n}.$ 
\begin{proof}
The cumulative distribution function is
$$F_{X_{1:n}}(x)= 1-\prod_{i=1}^{n}\left( 1-e^{-\theta x^{-\phi}}\right)^{\alpha_{i}}, \hspace{.25cm}x,\theta, \phi, \alpha_{i} >0.$$ The density functions  of $ X_{1:n} $ and  $ X^{*}_{1:n} ,$ respectively, are 
$$ f_{X_{1:n}}(x)=\theta \phi x^{-\phi-1} e^{-\theta x^{-\phi}}\left( \sum_{i=1}^{n}\alpha_{i}\right) \left(  1-e^{-\theta x^{-\phi}}\right)^{\sum_{i=1}^{n}\alpha_{i}-1}, $$
and $$ f_{X^{*}_{1:n}}(x)=\theta \phi x^{-\phi-1} e^{-\theta x^{-\phi}}\left( \sum_{i=1}^{n}\alpha^{*}_{i}\right) \left(  1-e^{-\theta x^{-\phi}}\right)^{\sum_{i=1}^{n}\alpha^{*}_{i}-1}, i=1,...,n. $$
Now, consider the function 
$$\displaystyle \dfrac{f_{X_{1:n}}(x)}{f_{X^{*}_{1:n}}(x)}=\dfrac{\sum_{i=1}^{n}\alpha_{i}}{\sum_{i=1}^{n}\alpha^{*}_{i}} \left(  1-e^{-\theta x^{-\phi}}\right)^{\sum_{i=1}^{n}\alpha_{i}-\sum_{i=1}^{n}\alpha^{*}_{i}}
$$
Taking logarithms and differentiating with respect to $x$ gives
$$\frac{d}{dx}\log\left(\frac{f_{X_{1:n}}(x)}{f_{X_{1:n}^*}(x)}\right)=\bigg(\sum \alpha_i-\sum \alpha_i^*\bigg)\frac{-\theta\phi x^{-\phi-1}e^{-\theta x^{-\phi}}}{1-e^{-\theta x^{-\phi}}}.$$
Thus, the derivative is nonnegative for all $x>0$ if and only if $\sum \alpha_i \le \sum \alpha_i^*$, completing the proof.
\end{proof}
\label{th9}
\end{theorem}

\begin{corollary}
 Let $ X_{1},\cdots.,X_{n} ,$ and $ X^{*}_{1},\cdots,X^{*}_{n} $  be two sets of non-negative independent random variables such that $ X_{i}\sim EG_{2}(\theta, \phi, \alpha_{i}),$ and  $ X^{*}_{i}\sim EG_{2}(\theta, \phi, \alpha_{i}^{*})$,
$ i=1,2,....,n$. If $\sum_{i=1}^{n}\alpha_{i}\leq \sum_{i=1}^{n}\alpha^{*}_{i} $ then
\begin{itemize}
    \item[(i)]  $ X^*_{1:n}\le_{mrl} X_{1:n},
\ \ \text{for all }x>0.$
\item[(ii)]   $ X^*_{1:n}\le_{icx} X_{1:n},
\ \ \text{for all }x>0.$
\end{itemize}
\end{corollary}

\begin{proof}
    The proof follows from the chain of implications existing between stochastic orders.
    \begin{itemize}
        \item[(i)]   Let $U$ and $V$ be two random variables. It is well-established in reliability theory (see Theorem 2.A.1; pp-83, \cite{shaked2007stochastic}) that  $U \leq_{lr} V \implies U \leq_{fr} V \implies U \leq_{mrl} V$. Now, in Theorem-\ref{th9}, we have established that $\sum_{i=1}^{n}\alpha_{i}\leq \sum_{i=1}^{n}\alpha^{*}_{i} \iff X^*_{1:n}\le_{lr} X_{1:n}$. By the aforementioned implication, this directly ensures $X^*_{1:n} \leq_{mrl} X_{1:n}$.
        
        \item[(ii)]  Similarly, it is known that $U \leq_{st} V \implies U \leq_{icx} V$ (see Theorem 4.A.55, p. 205 of \cite{shaked2007stochastic}). Since $X^*_{1:n} \leq_{lr} X_{1:n}$ implies $X^*_{1:n} \leq_{st} X_{1:n}$, and we have established (Theorem-\ref{th9}) that $\sum_{i=1}^{n}\alpha_{i}\leq \sum_{i=1}^{n}\alpha^{*}_{i} \iff X^*_{1:n}\le_{lr} X_{1:n}$, the condition $\sum \alpha_i \leq \sum \alpha^*_i$ is sufficient for the increasing convex ordering.
    \end{itemize}
\end{proof}

The failure rate function plays a central role in reliability modelling because it quantifies the instantaneous risk of system failure at any mission time. Since the lowest order statistic models a series system lifetime, it is important to know when two such systems satisfy a proportional failure structure. The following theorem addresses this issue for $EG_2$ components.

\begin{theorem}\label{thm:ph-minmax-alpha}
Let $ X_{1},\cdots.,X_{n} ,$ and $ X^{*}_{1},\cdots,X^{*}_{n} $  be two sets of non-negative independent random variables such that $ X_{i}\sim EG_{2}(\theta, \phi, \alpha_{i}),$ and  $ X^{*}_{i}\sim EG_{2}(\theta, \phi, \alpha_{i}^{*})$,
$ i=1,2,....,n$.  
Then the following statements hold:
\begin{enumerate}
\item[(i)] The lowest order statistics $X_{1:n}$ and $X^*_{1:n}$ satisfy a proportional failures relation.
\item[(ii)] In general, the highest order statistics $X_{n:n}$ and $X^*_{n:n}$ do not satisfy
a proportional failures relation.
\end{enumerate}
\end{theorem}

\begin{proof}
For each $i=1,\dots,n$, the reliability function and density of $X_i\sim EG_2(\theta,\phi, \alpha_i)$ are
\[
\overline{F}_{X_i}(x)=\Big(1-e^{-\theta x^{-\phi}}\Big)^{\alpha_i}, \qquad x>0,
\]
and
\[
f_{X_i}(x)=\alpha_i\theta\phi x^{-\phi-1}e^{-\theta x^{-\phi}}
\Big(1-e^{-\theta x^{-\phi}}\Big)^{\alpha_i-1}, \qquad x>0.
\]

(i) Since the variables are independent, the reliability function of the lowest order statistic is
\[
\overline{F}_{X_{1:n}}(x)=\prod_{i=1}^n \overline{F}_{X_i}(x)
=\prod_{i=1}^n \Big(1-e^{-\theta x^{-\phi}}\Big)^{\alpha_i}
=\Big(1-e^{-\theta x^{-\phi}}\Big)^{\sum_{i=1}^n \alpha_i}.
\]
Differentiating, we obtain
\[
f_{X_{1:n}}(x)
=
\sum_{i=1}^n \alpha_i~~ \Big(1-e^{-\theta x^{-\phi}}\Big)^{\sum_{i=1}^n \alpha_i-1}~~\,\Big(0-\theta\phi x^{-\phi-1}e^{-\theta x^{-\phi}}\Big)
.
\]
Therefore, the failure rate function of $X_{1:n}$ is
\begin{eqnarray}
    r_{X_{1:n}}(x)
&=&
\frac{f_{X_{1:n}}(x)}{\overline{F}_{X_{1:n}}(x)}\nonumber\\
&=&\frac{\sum_{i=1}^n \alpha_i~~ \Big(1-e^{-\theta x^{-\phi}}\Big)^{\sum_{i=1}^n \alpha_i-1}~~\,\Big(0-\theta\phi x^{-\phi-1}e^{-\theta x^{-\phi}}\Big)}{\Big(1-e^{-\theta x^{-\phi}}\Big)^{\sum_{i=1}^n \alpha_i}}\nonumber\\
&=&
\,\frac{-\theta\phi x^{-\phi-1}e^{-\theta x^{-\phi}}}
{1-e^{-\theta x^{-\phi}}} \sum_{i=1}^n \alpha_i.\nonumber
\end{eqnarray}

Similarly,
\[
r_{X^*_{1:n}}(x)
=
\,\frac{-\theta\phi x^{-\phi-1}e^{-\theta x^{-\phi}}}
{1-e^{-\theta x^{-\phi}}} \sum_{i=1}^n \alpha^*_i.
\]
Hence,
\[
\frac{r_{X_{1:n}}(x)}{r_{X^*_{1:n}}(x)}=\frac{\sum_{i=1}^n \alpha_i}{\sum_{i=1}^n \alpha_i^*}, \qquad x>0,
\]
which is independent of $x$. This proves that $X_{1:n}$ and $X^*_{1:n}$ are related through a
proportional failures model. 

(ii) Now consider the highest order statistic. Since the variables are independent,
\[
\overline{F}_{X_{n:n}}(x)
=
1-\prod_{i=1}^n F_{X_i}(x)
=
1-\prod_{i=1}^n \left[1-\Big(1-e^{-\theta x^{-\phi}}\Big)^{\alpha_i}\right].
\]
Consequently,
\[
f_{X_{n:n}}(x)
=
\sum_{i=1}^n
\alpha_i\theta\phi x^{-\phi-1}e^{-\theta x^{-\phi}}
\Big(1-e^{-\theta x^{-\phi}}\Big)^{\alpha_i-1}
\prod_{\substack{j=1\\j\ne i}}^n
\left[1-\Big(1-e^{-\theta x^{-\phi}}\Big)^{\alpha_j}\right].
\]
Hence the failure rate of $X_{n:n}$ is
\[
r_{X_{n:n}}(x)
=
\frac{
\displaystyle
\sum_{i=1}^n
\alpha_i\theta\phi x^{-\phi-1}e^{-\theta x^{-\phi}}
\Big(1-e^{-\theta x^{-\phi}}\Big)^{\alpha_i-1}
\prod_{\substack{j=1\\j\ne i}}^n
\left[1-\Big(1-e^{-\theta x^{-\phi}}\Big)^{\alpha_j}\right]
}{
\displaystyle
1-\prod_{i=1}^n \left[1-\Big(1-e^{-\theta x^{-\phi}}\Big)^{\alpha_i}\right]
}.
\]
An analogous representation holds for $r_{Y_{n:n}}(x)$. Unlike the lowest case, the above expression does not reduce to a common baseline failure
multiplied by a scalar depending only on $\sum_{i=1}^n\alpha_i$. Rather, it depends on the
entire vector $\boldsymbol{\alpha}=(\alpha_{1},\cdots\alpha_{n})$ through a nonlinear combination of products and sums.
Therefore, in general, the ratio $r_{X_{n:n}}(x)/r_{X^*_{n:n}}(x)$ is a nontrivial function of $x$
and is not constant on $(0,\infty)$. Thus, without further restrictions on
$\boldsymbol{\alpha}$ and $\boldsymbol{\alpha}^*$, no proportional failures relation can be asserted
for the highest order statistics. This completes the proof.
\end{proof}

\begin{remark}\label{rem:relative-aging-minmax}
The above theorem shows that, under common $(\theta,\phi)$, the lowest order statistic admits
a particularly simple relative aging interpretation with respect to the outer-shape parameters,
the effect of $\boldsymbol{\alpha}$ is entirely through the scalar $\sum_{i=1}^n\alpha_i$. Furthermore, the constant failure rate ratio ensures that the two systems are ordered according to \textit{relative aging orders}.
This phenomenon does not persist for the highest order statistic, whose failure rate depends on the
entire vector $\boldsymbol{\alpha}$ and therefore does not, in general, exhibit proportional failures. In the following, we have shown this with a counterexample.
\end{remark}

\begin{figure}[!h]
	\centering
	\includegraphics[scale=0.5]{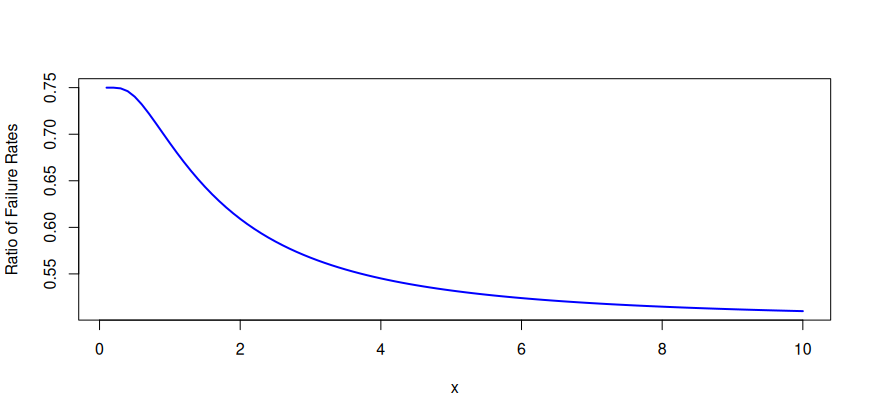}
\caption{Ratio of the failure rate functions $r_{X_{n:n}}(x)/r_{X^*_{n:n}}(x)$ for the $EG_2$ highest order statistics with $\boldsymbol{\alpha}=(1,3)$, $\boldsymbol{\alpha}^*=(2,2)$, and common $\theta = \phi = 1$. The non-constant trend confirms the absence of a proportional failures relationship for the highest  order statistics.}
    \label{xnn_ratio_failure}
\end{figure}

\begin{example}[Counterexample]
To illustrate the result in Theorem \ref{thm:ph-minmax-alpha}(ii), we consider a numerical case for $n=2$ with baseline parameters $\theta=1, \phi=1$. We compare two systems with different parameter vectors but identical sums:
\
where $\sum \alpha_i = \sum \alpha_i^* = 4$. Figure-\ref{xnn_ratio_failure} displays the ratio of the failure rate functions $r_{X_{n:n}}(x)/r_{X^*_{n:n}}(x)$. 

As observed, the ratio is strictly non-constant over the domain $x \in (0, 10)$, starting at approximately $0.75$ and decreasing toward $0.50$. This numerical evidence confirms that for the highest order statistic, the failure rate does not reduce to a simple proportional model. Instead, the aging behavior is sensitive to the individual components of the vector $\boldsymbol{\alpha}$, which further supports the need for more complex ordering tools for the highest order statistics.
\end{example}

In addition to stochastic dominance, reliability analysis often requires comparison of the spread, aging tendency, and variability of system lifetimes. The star order and dispersive order are useful tools for this purpose, since they provide deeper information than the usual stochastic order alone. For series systems, such results help in understanding not only which system lasts longer, but also which one exhibits greater uncertainty or stronger aging behaviour. The following theorem, therefore, strengthens the previous comparison by establishing star and dispersive ordering for the lowest order statistics.
\newpage
\begin{theorem}\label{thm:star-disp-min-alpha}
Let $ X_{1},\cdots.,X_{n} ,$ and $ X^{*}_{1},\cdots,X^{*}_{n} $  be two sets of non-negative independent random variables such that $ X_{i}\sim EG_{2}(\theta, \phi, \alpha_{i}),$ and  $ X^{*}_{i}\sim EG_{2}(\theta, \phi, \alpha_{i}^{*})$,
$ i=1,2,....,n$.  
Write
\[
\mathcal{A}=\sum_{i=1}^n \alpha_i,\qquad \mathcal{A}^*=\sum_{i=1}^n \alpha_i^*,  ~~~\text{ and set}\qquad c=\frac{\mathcal A}{\mathcal A^*}.
\]
If $\mathcal A\ge \mathcal A^*$, then
$X_{1:n}\le_* X^*_{1:n}
\qquad\text{and}\qquad
X_{1:n}\le_{disp} X^*_{1:n}.$
\end{theorem}

\begin{proof}
Since the components are independent and have common $(\theta,\phi)$, the distribution
functions of the lowest order statistics are
$$F_{X_{1:n}}(x)=1-\Big(1-e^{-\theta x^{-\phi}}\Big)^ \mathcal A,
\qquad
F_{X^*_{1:n}}(x)=1-\Big(1-e^{-\theta x^{-\phi}}\Big)^{ \mathcal A^*},
\qquad x>0.$$
To prove the star order, it is enough to show that $\frac{F_{X^*_{1:n}}^{-1}(F_{X_{1:n}}(x))}{x}$
is increasing in $x>0$. Let us consider the following

\begin{eqnarray}
h(x)&=& F_{X^*_{1:n}}^{-1}(F_{X_{1:n}}(x))\nonumber\\
F_{X^*_{1:n}}(h(x))&=& F_{X_{1:n}}(x) \nonumber\\
1-\Big(1-e^{-\theta h(x)^{-\phi}}\Big)^{\mathcal A^*}
&=&
1-\Big(1-e^{-\theta x^{-\phi}}\Big)^\mathcal A \nonumber\\
1-e^{-\theta h(x)^{-\phi}}
&=&
\Big(1-e^{-\theta x^{-\phi}}\Big)^{\mathcal A/ \mathcal A^*}\nonumber\\
&=&
\Big(1-e^{-\theta x^{-\phi}}\Big)^c \nonumber\\
h(x)
&=&
\left[
-\frac{1}{\theta}
\log\left(1-\Big(1-e^{-\theta x^{-\phi}}\Big)^c\right)
\right]^{-1/\phi}\nonumber\\
\frac{h(x)}{x}
&=&
\left[
\frac{\theta x^{-\phi}}
{-\log\left(1-\Big(1-e^{-\theta x^{-\phi}}\Big)^c\right)}
\right]^{1/\phi}.
\label{x1n_exp}
\end{eqnarray}
Now put $z=\theta x^{-\phi},\qquad
M_c(z)=-\log\left(1-\big(1-e^{-z}\big)^c\right),$
so that $\frac{h(x)}{x}
=
\left(\frac{z}{M_c(z)}\right)^{1/\phi}.$
Thus, it is sufficient to show that $z/M_c(z)$ is decreasing in $z>0$.
For this, note that $M_c(0)=0$, $c\ge 1$, and 
\[
M_c'(z)
=
\frac{c e^{-z}(1-e^{-z})^{\,c-1}}{1-(1-e^{-z})^c}=\gamma(c,1-e^{-z}) ~~~~~~~~~~~~~~\text{~~~~~~~~[from Lemma \ref{lemma4}]}.
\]
Now, using condition (3) of Lemma \ref{lemma4}, since $1-e^{-z}$ is increasing in $z$, it follows that $M_c'(z) = \gamma(c, 1-e^{-z})$ is increasing in $z$ for $c \geq 1$. Since $M_c(0)=0$, by the monotone form of L'Hospital's rule, it follows that $\frac{M_c(z)}{z}$ is increasing in $z > 0$. Consequently, $\frac{z}{M_c(z)}$ is decreasing in $z$, and since $z = \theta x^{-\phi}$ is a decreasing function of $x$, the ratio 
\[
\frac{h(x)}{x} = \left(\frac{z}{M_c(z)}\right)^{1/\phi}
\]
is increasing in $x$, which proves $X_{1:n} \le_* X^*_{1:n}$. 


Next, since $ \mathcal  A\ge \mathcal  A^*$, it can be shown similar to Theorem-\ref{th9} that $X_{1:n}\le_{lr} X^*_{1:n} \implies X_{1:n}\le_{st} X^*_{1:n}.$ Hence, using Theorem 4.B.2., p. 215 of \cite{shaked2007stochastic}, we conclude that $X_{1:n}\le_{disp} X^*_{1:n}.$
\end{proof}

For a parallel system, the highest order statistic determines the system lifetime, and hence any variability comparison for the highest has a direct reliability interpretation. In practical applications, homogeneous outer-shape settings arise naturally when the components are manufactured under similar operating environments or share common structural characteristics. Under this framework, star and dispersive orderings provide useful information on the variability and aging pattern of parallel-system lifetimes. The next theorem establishes these comparisons for the highest order statistic.

\begin{theorem}\label{thm:star-disp-max-alpha}
Let $ X_{1},\cdots.,X_{n} ,$ and $ X^{*}_{1},\cdots,X^{*}_{n} $  be two sets of non-negative independent random variables such that $ X_{i}\sim EG_{2}(\theta, \phi, \alpha),$ and  $ X^{*}_{i}\sim EG_{2}(\theta, \phi, \alpha^{*})$,
$ i=1,2,....,n$.  
If $\alpha\ge \alpha^*$, then
\[
X_{n:n}\le_* X^*_{n:n}
\qquad\text{and}\qquad
X_{n:n}\le_{disp} X^*_{n:n}.
\]
\end{theorem}

\begin{proof}
For $x>0$, the cumulative distribution functions of the highest order statistics are
\[
F_{X_{n:n}}(x)
=
\left[1-\Big(1-e^{-\theta x^{-\phi}}\Big)^{\alpha}\right]^n, \text{~~~~and~~~}  F_{X^*_{n:n}}(x)
=
\left[1-\Big(1-e^{-\theta x^{-\phi}}\Big)^{\alpha^*}\right]^n.
\]
Set $c=\frac{\alpha}{\alpha^*}\ge 1,
\quad
h(x):=F_{X^*_{n:n}}^{-1}(F_{X_{n:n}}(x)).$ Then we have the following

\begin{eqnarray}
    F_{X^*_{n:n}}(h(x))&=& F_{X_{n:n}}(x) \nonumber\\
    \left[1-\Big(1-e^{-\theta h(x)^{-\phi}}\Big)^{\alpha^*}\right]^n
&=&
\left[1-\Big(1-e^{-\theta x^{-\phi}}\Big)^{\alpha}\right]^n \nonumber\\
1-\Big(1-e^{-\theta h(x)^{-\phi}}\Big)^{\alpha^*}
&=&
1-\Big(1-e^{-\theta x^{-\phi}}\Big)^{\alpha} \nonumber\\
1-e^{-\theta h(x)^{-\phi}}
&=&
\Big(1-e^{-\theta x^{-\phi}}\Big)^{\alpha/\alpha^*} \nonumber\\
&=&
\Big(1-e^{-\theta x^{-\phi}}\Big)^c \nonumber\\
h(x)
&=&
\left[
-\frac{1}{\theta}
\log\left(1-\Big(1-e^{-\theta x^{-\phi}}\Big)^c\right)
\right]^{-1/\phi}\nonumber\\
\frac{h(x)}{x}
&=&
\left[
\frac{\theta x^{-\phi}}
{-\log\left(1-\Big(1-e^{-\theta x^{-\phi}}\Big)^c\right)}
\right]^{1/\phi}.
\end{eqnarray}

This is the same expression as in Equation-\ref{x1n_exp} in the proof of
Theorem~\ref{thm:star-disp-min-alpha}. Hence, by the same argument,
$h(x)/x$ is increasing in $x$, and therefore
\[
X_{n:n}\le_* X^*_{n:n}.
\]
Further, since $\alpha\ge \alpha^*$, we have $1-\Big(1-e^{-\theta x^{-\phi}}\Big)^{\alpha}
\ge
1-\Big(1-e^{-\theta x^{-\phi}}\Big)^{\alpha^*},$ and so we have 
\[
F_{X_{n:n}}(x)\ge F_{X^*_{n:n}}(x),\qquad x>0. \text{~~~~Thus~~~}~~~~~ X_{n:n}\le_{st} X^*_{n:n}.
\]
Hence, using Theorem 4.B.2., p. 215 of \cite{shaked2007stochastic}, we conclude that 
\[
X_{n:n}\le_{disp} X^*_{n:n}.
\]

\end{proof}

\begin{remark}
The conclusion of Theorem~\ref{thm:star-disp-max-alpha} is proved under the homogeneous
outer-shape setting $\alpha_1=\cdots=\alpha_n=\alpha$ and
$\alpha_1^*=\cdots=\alpha_n^*=\alpha^*$.
For heterogeneous outer-shape vectors
$\boldsymbol{\alpha}=(\alpha_1,\dots,\alpha_n)$ and
$\boldsymbol{\alpha}^*=(\alpha_1^*,\dots,\alpha_n^*)$,
the distribution of the highest order statistic does not reduce to a one-parameter form,
and a corresponding general proof is not immediate.
\end{remark}

The parameter $\phi$ governs an important aspect of the shape of the $EG_2$ distribution, and therefore changes in this parameter affect the failure-time pattern of the components. In reliability design, it is often of interest to understand how heterogeneity in such a shape parameter influences the performance of a \emph{series system}. Since the lifetime of a series system is the lowest order statistic, the next theorem studies the effect of majorization of the inner-shape parameter vector on the system lifetime in the sense of usual stochastic ordering. Note that, here, $\theta$ is the scale parameter, and $\alpha$ is the outer shape parameter.

\begin{theorem}
Let $ X_{1},\cdots.,X_{n} ,$ and $ X^{*}_{1},\cdots,X^{*}_{n} $  be two sets of non-negative independent random variables such that $ X_{i}\sim EG_{2}(\theta, \phi_{i}, \alpha),$ and  $ X^{*}_{i}\sim EG_{2}(\theta, \phi_{i}^{*}, \alpha)$,\\
$ i=1,2,....,n$. Then, if  $ (\phi_{1},\cdots,\phi_{n})\succeq^{m}(\phi^{*}_{1},\cdots,\phi^{*}_{n}) $, we have 
$ X_{1:n}\leq_{st} X_{1:n}^{*} .$

\begin{proof}
The cumulative distribution function and the reliability function of $ X_{1:n} ,$ respectively, are

$$F_{X_{1:n}}(x)= 1-\prod_{i=1}^{n}\left( 1-e^{-\theta x^{-\phi_{i}}}\right)^{\alpha},
$$and
$$\overline{F}_{X_{1:n}}(x)=\prod_{i=1}^{n}\left( 1-e^{-\theta x^{-\phi_{i}}}\right)^{\alpha}, \hspace{.25cm}x,\theta, \phi_{i}, \alpha >0, i=1,...,n.$$\\
To prove this theorem, it is enough to show the function, $\overline{F}_{X_{1:n}}(x)  $ is a\textit{ Schur-concave} function with respect to $ (\phi_{1},\cdots,\phi_{n}) $.
Now, the partial derivative of $\overline{F}_{X_{1:n}}(x) $ with respect to $ \phi_{i}, i=1,....,n. $ is \\

\begin{center}
$ \dfrac{\partial \overline{F}_{X_{1:n}}(x)}{\partial \phi_{i}}=-\theta \alpha(\log x)\overline{F}_{X_{1:n}}(x)  \dfrac{x^{-\phi_{i}}e^{-\theta x^{-\phi_{i}}}}{\left( 1-e^{-\theta x^{-\phi_{i}}}\right)}$\\
$\hspace{2cm} =  \alpha(\log x)\overline{F}_{X_{1:n}}(x)  \dfrac{e^{-\theta x^{-\phi_{i}}}\log e^{-\theta x^{-\phi_{i}}}}{\left( 1-e^{-\theta x^{-\phi_{i}}}\right)}$\\
$\hspace{.15cm} =\alpha(\log x)\overline{F}_{X_{1:n}}(x)g(\theta,\phi_{i}), $ 
\end{center}
where $ g(\theta,\phi_{i})=\dfrac{e^{-\theta x^{-\phi_{i}}}\log e^{-\theta x^{-\phi_{i}}}}{\left( 1-e^{-\theta x^{-\phi_{i}}}\right)},\hspace{.1cm} i=1,...,n.$   Choose $ u_{i}= e^{-\theta x^{-\phi_{i}}},$ it is clear that $  0<u_{i}<1 $ for $ i=1,...,n.$ Therefore,
using Lemma \ref{lemma3} we find that $ g(\theta,\phi_{i}) $ is decreasing in $ u_{i}=e^{-\theta x^{-\phi_{i}}} , i=1,...,n.$ Now consider, $u(\phi) = e^{-\theta x^{-\phi}}$, taking derivative with respect to $\phi$ we get $\displaystyle \frac{\partial u}{\partial \phi} = u \cdot \theta x^{-\phi} \log x.$
The sign of $\frac{\partial u}{\partial \phi}$ is determined by $\log x$, when $0 < x < 1$, $\log x < 0$ $\implies$ $u$ is decreasing in $\phi$, and when  $x > 1$, $\log x > 0$ $\implies$ $u$ increasing in $\phi$. Now let us consider two following cases:
\begin{enumerate}
 \item When $0 < x < 1$, then $u$ is decreasing in $\phi$, and since  $ g(\theta,\phi) $ is decreasing in $ u$ then using the chain rule of derivative we can get  $\frac{dg(\theta,\phi)}{d\phi}= \frac{\partial g(\theta,\phi)}{\partial u}\times \frac{\partial u}{d\phi}>0 $, hence $g(\theta,\phi)$ is increasing in $ \phi$.
    \item  Similarly, for $x > 1$,  $u$ is increasing in $\phi$ and again since  $ g(\theta,\phi) $ is decreasing in $ u$ then using the chain rule of derivative we can get   $\frac{dg(\theta,\phi)}{d\phi}= \frac{\partial g(\theta,\phi)}{\partial u}\times \frac{\partial u}{d\phi}<0 $, hence $g(\theta,\phi)$ is decreasing in $ \phi$.
\end{enumerate}

Now, without loss of generality, let us consider $
\phi_{i}\geq\phi_{j}, $ and using the immediate above conditions we have
$$\phi_{i}-\phi_{j})\left( \dfrac{\partial \overline{F}_{X_{1:n}}(x)}{\partial \phi_{i}}- \dfrac{\partial \overline{F}_{X_{1:n}}(x)}{\partial \phi_{j}}\right)=\begin{cases}
  (\phi_{i}-\phi_{j})\alpha(\log x)\overline{F}_{X_{1:n}}(x)  \left[ g(\theta,\phi_{i})- g(\theta,\phi_{j})\right]\leq0 . & \text{if } x \in (0,1] \\
 (\phi_{i}-\phi_{j})\alpha(\log x)\overline{F}_{X_{1:n}}(x)  \left[ g(\theta,\phi_{i})- g(\theta,\phi_{j})\right]<0 . & \text{if } x > 1
\end{cases}$$
Therefore for fixed $ x>0 ,$ we get 
$(\phi_{i}-\phi_{j})\left( \dfrac{\partial \overline{F}_{X_{1:n}}(x)}{\partial \phi_{i}}- \dfrac{\partial \overline{F}_{X_{1:n}}(x)}{\partial \phi_{j}}\right) \leq0 . $ Using Theorem \ref{th1}, we say that  $\overline{F}_{X_{1:n}}(x) $ is a \textit{Schur-concave} function with respect to  $ (\phi_{1},\cdots,\phi_{n}).$ So, we obtain $\overline{F}_{X^{*}_{1:n}}(x)\geq\overline{F}_{X_{1:n}}(x)  $ i.e.,  $ X_{1:n}\leq_{st} X_{1:n}^{*},$ as required.
\end{proof}
\label{th10}
\end{theorem}

\begin{corollary}
   Let $ X_{1},\cdots.,X_{n} ,$ and $ X^{*}_{1},\cdots,X^{*}_{n} $  be two sets of non-negative independent random variables such that $ X_{i}\sim EG_{2}(\theta, \phi_{i}, \alpha),$ and  $ X^{*}_{i}\sim EG_{2}(\theta, \phi_{i}^{*}, \alpha)$,
$ i=1,2,....,n$. Then, if  $ (\phi_{1},\cdots,\phi_{n})\succeq^{m}(\phi^{*}_{1},\cdots,\phi^{*}_{n}) $, we have 
$ X_{1:n}\leq_{icx} X_{1:n}^{*} .$
\end{corollary}

\begin{proof}
        The proof is immediate from the 
        Theorem-4.A.55 at pp-205 of \cite{shaked2007stochastic}.
\end{proof}

Unlike the series-system, the effect of the inner-shape parameter on the lifetime of a \emph{parallel system} is more subtle. Since the parallel-system lifetime depends on the highest order statistic, a general stochastic ordering need not always hold under parameter majorization. From a reliability point of view, it is important to identify such limitations because they indicate that stronger assumptions are required before one system can be declared superior to another. The following proposition highlights this non-comparability phenomenon.

\begin{proposition}
The highest order statistics are not comparable in the usual stochastic ordering with respect to the inner shape parameter $ \phi. $  
\end{proposition}
\begin{proof}
Let $X_i \sim EG_2(\theta,\phi_i,\alpha)$, $i=1,\dots,n$, be independent random variables.
The cumulative distribution function of $X_i$ is
\[
F_i(x)
=
1-\left(1-e^{-\theta x^{-\phi_i}}\right)^{\alpha},
\quad x>0.
\]

The reliability function of the highest order statistic $X_{n:n}$ is
\[
\overline F_{X_{n:n}}(x)
=
P(X_{n:n}>x)
=
1-\prod_{i=1}^n F_i(x)
=
1-
\prod_{i=1}^n
\left[
1-\left(1-e^{-\theta x^{-\phi_i}}\right)^{\alpha}
\right].
\]

For notational convenience, define $H_i(x)
=
1-\left(1-e^{-\theta x^{-\phi_i}}\right)^{\alpha},
\quad 0<H_i(x)<1.$

Then $\overline F_{X_{n:n}}(x)
=
1-\prod_{i=1}^n H_i(x).$
Differentiating $\overline F_{X_{n:n}}(x)$ with respect to $\phi_i$, we obtain
\begin{align*}
    \frac{\partial \overline F_{X_{n:n}}(x)}{\partial \phi_i}&=-
\left(\prod_{k\neq i} H_k(x)\right) \frac{\partial H_i(x)}{\partial \phi_i}\\
&=-\alpha(1-e^{-\theta x^{-\phi_i}})^{\alpha-1}
\frac{\partial}{\partial \phi_i}
(1-e^{-\theta x^{-\phi_i}})\prod_{k\neq i} H_k(x) \\
&=- \alpha\theta x^{-\phi_i}\log x\,
e^{-\theta x^{-\phi_i}}
(1-e^{-\theta x^{-\phi_i}})^{\alpha-1}
\prod_{k\neq i} H_k(x).
\end{align*}


All factors in the above expression are strictly positive except $\log x$.
Consequently,
\[
\operatorname{sign}
\left(
\frac{\partial \overline F_{X_{n:n}}(x)}{\partial \phi_i}
\right)
=
\operatorname{sign}(\log x).
\]

That is,
\[
\frac{\partial \overline F_{X_{n:n}}(x)}{\partial \phi_i}
\begin{cases}
>0, & 0<x<1,\\[4pt]
<0, & x>1.
\end{cases}
\]
Since the sign of the partial derivatives changes at $x=1$, the reliability
function $\overline F_{X_{n:n}}(x)$ is neither Schur-convex nor Schur-concave in
$(\phi_1,\dots,\phi_n)$ for all $x>0$.
Hence, majorization of the inner shape parameter vector does not imply
a usual stochastic order between $X_{n:n}$ and $X_{n:n}^*$ in general.
\end{proof}

\begin{figure}[!h]
	\centering
	\includegraphics[scale=0.8]{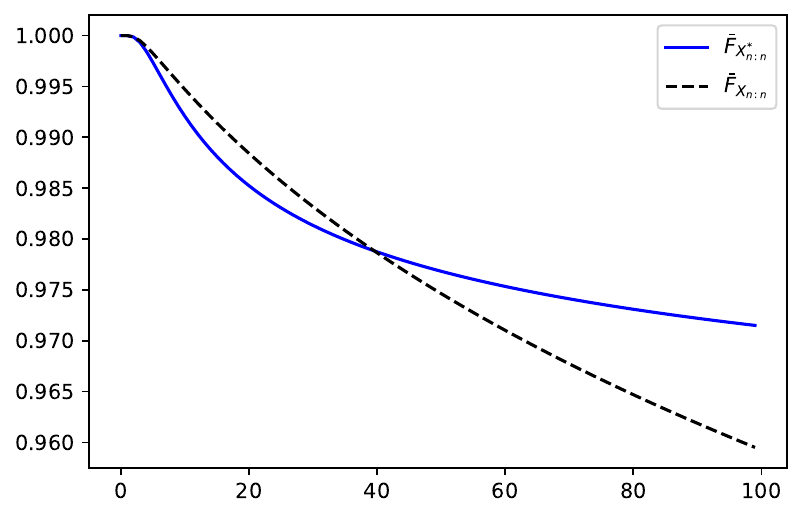}
	\caption{Presents the plot of $ \overline{F}_{X_{n:n}} $ and $ \overline{F}_{X^{*}_{n:n}}. $ }\label{A}
\end{figure}

In the next example, we have shown that we can not compare the highest order statistics in the usual stochastic ordering with respect to inner shape parameter $ \phi. $
\begin{example}
Let $ X_{1},X_{2},X_{3} $ and $ X^{*}_{1},X^{*}_{2},X^{*}_{3} $  be two sets of non-negative independent random variables such that $ X_{i}\sim EG_{2}(5, \phi_{i},2),$ and  $ X^{*}_{i}\sim EG_{2}(5, \phi_{i}^{*}, 2)$, $ i=1,2,3.$ Let us choose $ \phi_{1}=0.1, \phi_{2}=1.14,  \phi_{3}=0.3 $ and $ \phi^{*}_{1}=0.6, \phi^{*}_{2}=0.9,  \phi^{*}_{3}=0.04.$ It is clear that $ (\phi_{1},\cdots\phi_{n})\succeq^{m}(\phi^{*}_{1},\cdots\phi^{*}_{n}) .$ Now, from the Figure \ref{A}, we observe that the the graph of reliability functions of $\overline{F}_{X^{*}_{1:n}}(x)$ and $\overline{F}_{X_{1:n}}(x)  $cross each other. This implies that  $ X_{n:n} \ngeq_{st} X^{*}_{n:n}.$
\end{example}

\subsection{ Results based on multivariate chain majorization}
 \label{subsec2}
 Let $ X_{1},\cdots,X_{n} $ be a set of independent random variables and each $ X_{i} $ follows $ EG_{2} $ distributions i.e., $ X_{i}\sim EG_{2}(\theta_{i}, \phi_{i}, \alpha_{i}), i=1,2,....,n$. Furthermore, let $ X^{*}_{1},\cdots,X^{*}_{n} $ be another set of independent random variables such that $ X^{*}_{i}\sim EG_{2}(\theta_{i}^{*}, \phi_{i}^{*}, \alpha_{i}^{*}),$ $ i=1,2,....,n$. Assume $\phi_{1}=\cdots\phi_{n}=\phi^{*}_{1}=\cdots=\phi^{*}_{n}$. Here we deal with two systems (parallel, series) with independent $ EG_{2} $ components having different scale and outer shape parameters, i.e., different $\theta_{i} $'s, $ \alpha_{i}$'s for $ i=1,2,....,n. $\\

In many engineering applications, component parameters may vary simultaneously in more than one direction, for example, through changes in both \textit{scale} and \textit{outer shape} parameters with fixed \textit{inner shape parameter}, $\phi.$ In such settings, \emph{multivariate chain majorization} provides a natural framework for comparing heterogeneous systems. Since the lowest order statistic gives the lifetime of a series system, the following theorem is important because it extends the earlier univariate comparisons to a more flexible matrix-parameter setting, thereby allowing broader reliability comparisons between complex systems.
 
 \begin{theorem}
Let $ X_{1}, X_{2} $ and  $ X^{*}_{1}, X^{*}_{2} $ be two pairs of non-negative independent random variables such that $ X_{i}\sim EG_{2}(\theta_{i}, \phi, \alpha_{i}),$ and  $ X^{*}_{i}\sim EG_{2}(\theta_{i}^{*}, \phi, \alpha_{i}^{*})$ for $ i=1,2$. Then, 

 \(
 \begin{bmatrix}
 \alpha_{1}       & \alpha_{2}  \\
 \theta_{1}       & \theta_{2} \\
 
 \end{bmatrix}
 \) $ \in \mathcal{S}_{2} $ we have
 \(
 \begin{bmatrix}
 \alpha_{1}       & \alpha_{2}  \\
 \theta_{1}       & \theta_{2} \\
 
 \end{bmatrix}
 \)$ \gg $ \(
 \begin{bmatrix}
 \alpha^{*}_{1}       & \alpha^{*}_{2}  \\
 \theta^{*}_{1}       & \theta^{*}_{2} \\
 
 \end{bmatrix}
 \)$ \implies X^{*}_{1:2}\geq_{st}X_{1:2} .$

 \begin{proof}
 The reliability function of $ X_{1:2} $ is given by, for $ x>0 $
 $$\overline{F}_{X_{1:2}}(x)=\prod_{i=1}^{2}\left( 1-e^{-\theta_{i} x^{-\phi}}\right)^{\alpha_{i}},  \theta_{i}, \phi, \alpha_{i}>0,\hspace{.2cm} i=1,2 .$$
 It is easy to prove that the function
 $\overline{F}_{X_{1:2}}(x) $ is permutation invariant with respect to $ (\theta_{i},\alpha_{i}), i=1,2.$  So, the condition (1) of Theorem \ref{th2} is satisfied. 
 Next, we need to show that condition (2) of Theorem \ref{th2} also satisfies. Now, consider $\overline{\theta}=(\theta_{1},\theta_{2})  $ and $ \overline{\alpha}=(\alpha_{1},\alpha_{2}) $
 and let us define a function
 
 \begin{equation}
 G(\overline{\theta},\overline\alpha)=G_{1}(\overline{\theta},\overline\alpha)+G_{2}(\overline{\theta},\overline\alpha)
 \end{equation}
 Where, $$G_{1}(\overline{\theta},\overline\alpha)=(\alpha_{1}-\alpha_{2})\left( 
 \dfrac{\partial \overline{F}_{X_{1:2}}(x)}{\partial \alpha_{1}}-\dfrac{\partial \overline{F}_{X_{1:2}}(x)}{\partial \alpha_{2}} \right), $$
 $$G_{2}(\overline{\theta},\overline\alpha)=(\theta_{1}-\theta_{2})\left( 
 \dfrac{\partial \overline{F}_{X_{1:2}}(x)}{\partial \theta_{1}}-\dfrac{\partial \overline{F}_{X_{1:2}}(x)}{\partial \theta_{2}} \right). $$ The partial derivative of $ \overline{F}_{X_{1:2}}(x) $ with respect to  $ \overline\alpha $ gives
 $$ \dfrac{\partial \overline{F}_{X_{1:2}}(x)}{\partial \alpha_{i}}=\overline{F}_{X_{1:2}}(x) \log  \left( 1-e^{-\theta_{i} x^{-\phi}}\right), \hspace{.3cm}i=1,2. $$ Therefore 
 $$G_{1}(\overline{\theta},\overline\alpha)=(\alpha_{1}-\alpha_{2})\overline{F}_{X_{1:2}}(x) \left[ \log  \left( 1-e^{-\theta_{1} x^{-\phi}}\right)-\log  \left( 1-e^{-\theta_{2} x^{-\phi}}\right)\right] .$$ Our assumption is $ (\overline{\theta},\overline\alpha)\in \mathcal{S}_{2} .$ So, we have $ (\alpha_{1}-\alpha_{2})(\theta_{1}-\theta_{2})\leq0 $, this implies that either $\alpha_{1}\geq\alpha_{2},$ $\theta_{1}\leq\theta_{2} $ or, $\alpha_{1}\leq\alpha_{2},$ $\theta_{1}\geq\theta_{2} $. We choose the case when, $\alpha_{1}\leq\alpha_{2},$ $\theta_{1}\geq\theta_{2} $ for our proof. For the case when, $\alpha_{1}\geq\alpha_{2},$ $\theta_{1}\leq\theta_{2} $ the proof is quite similar. Now, we know that $ \log  \left( 1-e^{-\theta x^{-\phi}}\right) $ is increasing with respect to $ \theta .$ Then, we have  $$ \log  \left( 1-e^{-\theta_{1} x^{-\phi}}\right)\geq   \log  \left( 1-e^{-\theta_{2} x^{-\phi}}\right)  $$ and by assumption $\alpha_{1}\leq\alpha_{2}.$ This implies that $ G_{1}(\overline{\theta},\overline\alpha)\leq 0 .$\\
 
 The partial derivative of $ \overline{F}_{X_{1:2}}(x) $ with respect to  $ \overline\theta $ gives
 $$ \dfrac{\partial \overline{F}_{X_{1:2}}(x)}{\partial \theta_{i}}=x^{-\phi}\overline{F}_{X_{1:2}}(x)
 \dfrac{\alpha_{i}e^{-\theta_{i} x^{-\phi}}}{\left( 1-e^{-\theta_{i} x^{-\phi}}\right)},\hspace{.3cm}i=1,2.$$ Therefore 
 $$
 G_{2}(\overline{\theta},\overline\alpha)=(\theta_{1}-\theta_{2})x^{-\phi}\overline{F}_{X_{1:2}}(x)\left[ 
 \dfrac{\alpha_{1}e^{-\theta_{1} x^{-\phi}}}{\left( 1-e^{-\theta_{1} x^{-\phi}}\right)}-\dfrac{\alpha_{2}e^{-\theta_{2} x^{-\phi}}}{\left( 1-e^{-\theta_{2} x^{-\phi}}\right)}\right] .
 $$Our assumption is $\alpha_{1}\leq\alpha_{2},$
 $\theta_{1}\geq\theta_{2} .$ Now, the function $ \dfrac{\alpha e^{-\theta x^{-\phi}}}{\left( 1-e^{-\theta x^{-\phi}}\right)} $ is decreasing function with respect to $ \theta .$ Therefore we have $$\dfrac{\alpha_{1}e^{-\theta_{1} x^{-\phi}}}{\left( 1-e^{-\theta_{1} x^{-\phi}}\right)}\leq \dfrac{\alpha_{2}e^{-\theta_{2} x^{-\phi}}}{\left( 1-e^{-\theta_{2} x^{-\phi}}\right)}.$$ This implies that 
 $ G_{2}(\overline{\theta},\overline\alpha)\leq 0 .$ Since  $G_{1}(\overline{\theta},\overline\alpha)\leq 0 $ and $ G_{2}(\overline{\theta},\overline\alpha)\leq 0 $. Using (2) we conclude that $ G(\overline{\theta},\overline\alpha)\leq 0 ,$ and so the condition (2) of Theorem \ref{th2} is satisfied. Hence, the theorem follows.
 \end{proof}
 
 \label{th11}
 \end{theorem}

Parallel systems often exhibit more complicated reliability behavior than series systems, especially when several model parameters vary jointly. Therefore, obtaining stochastic comparisons of the highest order statistic under multivariate majorization is useful for understanding how joint parameter heterogeneity influences the lifetime of a parallel system. The following theorem addresses this problem and provides a reliability interpretation of the parameter matrix ordering in terms of system performance. Here, we fix the \textit{inner shape parameter}, $\phi.$
 \begin{theorem}
 Let $ X_{1}, X_{2} $ and  $ X^{*}_{1}, X^{*}_{2} $ be two pairs of non-negative independent random variables such that $ X_{i}\sim EG_{2}(\theta_{i}, \phi, \alpha_{i}),$ and  $ X^{*}_{i}\sim EG_{2}(\theta_{i}^{*}, \phi, \alpha_{i}^{*})$ for $ i=1,2$. Then, if
 \(
 \begin{bmatrix}
 \alpha_{1}       & \alpha_{2}  \\
 \theta_{1}       & \theta_{2} \\
 
 \end{bmatrix}
 \) $ \in \mathcal{T}_{2} $ we have
 \(
 \begin{bmatrix}
 \alpha_{1}       & \alpha_{2}  \\
 \theta_{1}       & \theta_{2} \\
 
 \end{bmatrix}
 \)$ \gg $ \(
 \begin{bmatrix}
 \alpha^{*}_{1}       & \alpha^{*}_{2}  \\
 \theta^{*}_{1}       & \theta^{*}_{2} \\
 
 \end{bmatrix}
 \)$ \implies X_{2:2}\geq_{st}X^{*}_{2:2} .$

 \begin{proof}
 The cumulative distribution function of $ X_{2:2} $ is given by, for $ x>0 $
 $$F_{X_{2:2}}(x)= \prod_{i=1}^{2}\left[ 1-\left( 1-e^{-\theta_{i} x^{-\phi}}\right)^{\alpha_{i}}\right], \theta_{i}, \phi, \alpha_{i}>0,\hspace{.2cm} i=1,2 .$$
 
 It is easy to check for fixed $ x>0, $  ${F}_{X_{2:2}}(x) $ is permutation invariant with respect to $ (\theta_{i},\alpha_{i}), i=1,2 .$ So, the condition (1) of Theorem \ref{th2} is satisfied. 
 Next, our claim is That condition (2) of Theorem  \ref{th2} also satisfies.  Now, consider $\overline\theta=(\theta_{1},\theta_{2})  $ and $ \overline\alpha=(\alpha_{1},\alpha_{2}) $ and 
 let us define a function
 
 \begin{equation}
 H(\overline\theta,\overline\alpha)=H_{1}(\overline\theta,\overline\alpha)+H_{2}(\overline\theta,\overline\alpha)
 \end{equation}
 Where, $$H_{1}(\overline\theta,\overline\alpha)=(\alpha_{1}-\alpha_{2})\left( 
 \dfrac{\partial F_{X_{2:2}}(x)}{\partial \alpha_{1}}-\dfrac{\partial F_{X_{2:2}}(x)}{\partial \alpha_{2}} \right), $$
 $$H_{2}(\overline\theta,\overline\alpha)=(\theta_{1}-\theta_{2})\left( 
 \dfrac{\partial F_{X_{2:2}}(x)}{\partial \theta_{1}}-\dfrac{\partial F_{X_{2:2}}(x)}{\partial \theta_{2}} \right). $$
 The partial derivative of $ F_{X_{2:2}}(x) $ with respect to  $\overline \alpha $ gives
 $$ \dfrac{\partial F_{X_{2:2}}(x)}{\partial \alpha_{i}}=- F_{X_{2:2}}(x)\dfrac{\left( 1-e^{-\theta_{i} x^{-\phi}}\right)^{\alpha_{i}}\log  \left( 1-e^{-\theta_{i} x^{-\phi}}\right) }{1-\left( 1-e^{-\theta_{i} x^{-\phi}}\right)^{\alpha_{i}}}, \hspace{.2cm}i=1,2 .$$ Let $ u_{i}= \left( 1-e^{-\theta_{i} x^{-\phi}}\right), i=1,2.$ Therefore we have
 $$ \dfrac{\partial F_{X_{2:2}}(x)}{\partial \alpha_{i}}=- F_{X_{2:2}}(x)\dfrac{u_{i}^{\alpha_{i}}\log  u_{i} }{1-u_{i}^{\alpha_{i}}}=- F_{X_{2:2}}(x)\eta(\alpha_{i},u_{i}),$$ where $ \eta(\alpha,u)=\dfrac{u^{\alpha}\log  u }{1-u^{\alpha}} .$ From Lemma \ref{lemma3},  it can be shown that  $ \eta(\alpha,u)=\eta(\alpha,1-e^{-\theta x^{-\phi}}) $ is increasing in $ \alpha $ for fixed $ \theta ,$ and is decreasing in $ \theta $ for fixed $ \alpha .$
 Now, by assumption $( \overline\theta,\overline\alpha) \in \mathcal{T}_{2}.$  So, we have $ (\alpha_{1}-\alpha_{2})(\theta_{1}-\theta_{2})\leq0 $, and $\alpha_{1},\alpha_{2}\geq$1, this implies that either $\alpha_{1}\geq\alpha_{2},$ $\theta_{1}\leq\theta_{2} $ or, $\alpha_{1}\leq\alpha_{2},$ $\theta_{1}\geq\theta_{2} $. We choose the case when, $\alpha_{1}\geq\alpha_{2}\geq 1,$ $\theta_{1}\leq\theta_{2} $ for our proof. For the case when, $\alpha_{2}\geq\alpha_{1}\geq 1,$ $\theta_{2}\leq\theta_{1} $ the proof is quite similar. 

 Therefore, we can conclude that
 $$H_{1}(\overline\theta,\overline\alpha)=(\alpha_{1}-\alpha_{2}) F_{X_{2:2}}(x)\left[  \eta(\alpha_{2},1-e^{-\theta_{2} x^{-\phi}})-\eta(\alpha_{1},1-e^{-\theta_{1} x^{-\phi}})     \right] \leq0.$$
 
 On the other hand the partial derivative of $ F_{X_{2:2}}(x) $ with respect to  $ \overline\theta $ gives
 $$ \dfrac{\partial F_{X_{2:2}}(x)}{\partial \theta_{i}}=-x^{-\phi}F_{X_{2:2}}(x)
 \dfrac{\alpha_{i} e^{-\theta_{i} x^{-\phi}} \left(  1-e^{-\theta_{i} x^{-\phi}}\right)^{\alpha_{i}-1}}{1-\left( 1-e^{-\theta_{i} x^{-\phi}}\right)^{\alpha_{i}}},\hspace{.2cm}i=1,2.$$
 Since $ u_{i}= \left( 1-e^{-\theta_{i} x^{-\phi}}\right), i=1,2.$ therefore we have 
 $$ \dfrac{\partial F_{X_{2:2}}(x)}{\partial \theta_{i}}=-x^{-\phi}F_{X_{2:2}}(x)
 \dfrac{\alpha_{i} (1-u_{i}) u_{i}^{\alpha_{i}-1}}{1-u_{i}^{\alpha_{i}}},\hspace{.2cm} i=1,2.$$

 For $\alpha_{1}\geq\alpha_{2}\geq 1,$ $\theta_{1}\leq\theta_{2} $ we have 
 
 $$
 H_{2}(\overline\theta,\overline\alpha)=(\theta_{1}-\theta_{2})x^{-\phi}F_{X_{2:2}}(x)\left[  \gamma(\alpha_{2},1-e^{-\theta_{2} x^{-\phi}})- \gamma(\alpha_{1},1-e^{-\theta_{1} x^{-\phi}}) \right]  \leq0
 ,$$ since $\gamma(\alpha,1-e^{-\theta x^{-\phi}})$ is increasing in $ \theta $ for fixed $ \alpha\geq1 ,$ and  is decreasing in $ \alpha$ for fixed $ \theta ,$ where $ \gamma(\alpha,u)=\dfrac{\alpha (1-u) u^{\alpha-1}}{1-u^{\alpha}},$ from Lemma \ref{lemma4}.  Therefore $ H_{1}(\overline\theta,\overline\alpha)\leq 0 $ and $ H_{2}(\overline\theta,\overline\alpha)\leq 0 $. So, from (3), we conclude that $ H(\overline\theta,\overline\alpha)\leq 0 .$ So, our claim for the condition (2) of Theorem \ref{th2} is satisfied. This completes the proof of the theorem.
 
 \end{proof}
 \label{th12}
 \end{theorem}

  Now, in the following, we construct two examples to verify our results from Theorem \ref{th11} and Theorem \ref{th12}.

 \begin{example}
 Let $ X_{1}, X_{2} $ be a pair of independent random variables such that $ X_{i}\sim EG_{2}(\theta_{i}, \phi, \alpha_{i}),\hspace{.2cm}  i=1,2$. Furthermore, let $ X^{*}_{1}, X^{*}_{2} $ be an another pair of independent random variables with $ X^{*}_{i}\sim EG_{2}(\theta_{i}^{*}, \phi, \alpha_{i}^{*}),\hspace{.2cm}i=1,2.$ Set\\
 
 \(
 \begin{bmatrix}
 \alpha_{1}       & \alpha_{2}  \\
 \theta_{1}       & \theta_{2} \\
 
 \end{bmatrix}
 \)$ = $ \(
 \begin{bmatrix}
 0.54       & 0.66  \\
 1.7       &1.4 \\
 
 \end{bmatrix}
 \) and  \(
 \begin{bmatrix}
 \alpha^{*}_{1}       & \alpha^{*}_{2}  \\
 \theta^{*}_{1}       & \theta^{*}_{2} \\
 
 \end{bmatrix}
 \)$ = $ \(
 \begin{bmatrix}
 0.5       &0 .7  \\
 1.8       &1.3 \\
 
 \end{bmatrix}
 \).\\

 We can see that both the matrices  belong in $ \mathcal{S}_{2}.$ Consider a T-transform matrix:
 \begin{center}
 $ T_{0.8}=0.8$\(
 \begin{bmatrix}
 1     & 0  \\
 0       &1 \\
 
 \end{bmatrix}
 \)$ + 0.2$ \(
 \begin{bmatrix}
 0     & 1  \\
 1       &0 \\
 
 \end{bmatrix}
 \)$ = $
 \(
 \begin{bmatrix}
 0.8      & 0.2  \\
 0.2       &0.8 \\
 
 \end{bmatrix}
 \), \end{center}
 
 and we have \(
 \begin{bmatrix}
 \alpha_{1}       & \alpha_{2}  \\
 \theta_{1}       & \theta_{2} \\
 
 \end{bmatrix}
 \)$ = $ \(
 \begin{bmatrix}
 \alpha^{*}_{1}       & \alpha^{*}_{2}  \\
 \theta^{*}_{1}       & \theta^{*}_{2} \\
 
 \end{bmatrix}
 \)$ T_{0.8} ,$ which satisfies the Definition \ref{df4},  so we have
 \begin{center}

 \(
 \begin{bmatrix}
 \alpha_{1}       & \alpha_{2}  \\
 \theta_{1}       & \theta_{2} \\
 
 \end{bmatrix}
 \)$ \gg $ \(
 \begin{bmatrix}
 \alpha^{*}_{1}       & \alpha^{*}_{2}  \\
 \theta^{*}_{1}       & \theta^{*}_{2} \\
 
 \end{bmatrix}
 \).
 \end{center}
 Therefore, from Theorem \ref{th11}, we get $  X_{1:2}\leq_{st}X^{*}_{1:2},$ for any $ x>0. $
 \end{example}
 
 \begin{example}
Let $ X_{1}, X_{2} $ be a pair of independent random variables such that $ X_{i}\sim EG_{2}(\theta_{i}, \phi, \alpha_{i}), \hspace{.2cm} i=1,2$. Furthermore, let $ X^{*}_{1}, X^{*}_{2} $ be an another pair of independent random variables with $ X^{*}_{i}\sim EG_{2}(\theta_{i}^{*}, \phi, \alpha_{i}^{*}),\hspace{.2cm} i=1,2.$ Set\\
 
 \(
 \begin{bmatrix}
 \alpha_{1}       & \alpha_{2}  \\
 \theta_{1}       & \theta_{2} \\
 
 \end{bmatrix}
 \)$ = $ \(
 \begin{bmatrix}
 2.34       & 2.26  \\
 1.32       &1.38 \\
 
 \end{bmatrix}
 \) and  \(
 \begin{bmatrix}
 \alpha^{*}_{1}       & \alpha^{*}_{2}  \\
 \theta^{*}_{1}       & \theta^{*}_{2} \\
 
 \end{bmatrix}
 \)$ = $ \(
 \begin{bmatrix}
 2.1       & 2.5  \\
 1.5       &1.2 \\
 
 \end{bmatrix}
 \).\\

 We can see that both the matrices  belong in $ \mathcal{T}_{2}.$
 Consider a T-transform matrix:
 \begin{center}

 $ T_{0.6}=0.4$\(
 \begin{bmatrix}
 1     & 0  \\
 0       &1 \\
 
 \end{bmatrix}
 \)$ + 0.6$ \(
 \begin{bmatrix}
 0     & 1  \\
 1       &0 \\
 
 \end{bmatrix}
 \)$ = $
 \(
 \begin{bmatrix}
 0.4      & 0.6  \\
 0.6       &0.4 \\
 
 \end{bmatrix}
 \), 
 \end{center}
 
 and we have \(
 \begin{bmatrix}
 \alpha_{1}       & \alpha_{2}  \\
 \theta_{1}       & \theta_{2} \\
 
 \end{bmatrix}
 \)$ = $ \(
 \begin{bmatrix}
 \alpha^{*}_{1}       & \alpha^{*}_{2}  \\
 \theta^{*}_{1}       & \theta^{*}_{2} \\
 
 \end{bmatrix}
 \)$ T_{0.4} .$ \\
 
 So, according to Definition \ref{df4} we have  
 \begin{center}

 \(
 \begin{bmatrix}
 \alpha_{1}       & \alpha_{2}  \\
 \theta_{1}       & \theta_{2} \\
 
 \end{bmatrix}
 \)$ \gg $ \(
 \begin{bmatrix}
 \alpha^{*}_{1}       & \alpha^{*}_{2}  \\
 \theta^{*}_{1}       & \theta^{*}_{2} \\
 
 \end{bmatrix}
 \).\\
 \end{center}
 
 Therefore, from Theorem \ref{th12}, we get $ X_{2:2}\geq_{st}X^{*}_{2:2},$ for any $ x>0. $
 \end{example}

 One-step majorization results are useful, but in practice, system comparisons often involve a sequence of transformations of the parameter matrix rather than a single one. Thus, from a reliability perspective, it is important to know whether stochastic comparisons remain valid under repeated transformations. The next theorems extend the previous result to more general ( $ n>2$) chain-majorization settings, thereby making the comparison framework applicable to a wider class of series-system configurations.
   
 \begin{theorem}
  Let $ X_{1},...., X_{n} $ and  $ X^{*}_{1},...., X^{*}_{n} $ be two sets of non-negative independent random variables such that $ X_{i}\sim EG_{2}(\theta_{i}, \phi, \alpha_{i}),$ and  $ X^{*}_{i}\sim EG_{2}(\theta_{i}^{*}, \phi, \alpha_{i}^{*})$ for $ i=1,...,n$.
 \begin{enumerate}
 \item  Assume that \(
 \begin{bmatrix}
 \alpha_{1}\hspace{.3cm}....&\alpha_{n}  \\
 \theta_{1} \hspace{.3cm}....& \theta_{n} \\
 
 \end{bmatrix}
 \) $ \in \mathcal{T}_{n} ,$ and
 
 \(
 \begin{bmatrix}
 \alpha^{*}_{1}\hspace{.3cm}....& \alpha^{*}_{n}  \\
 \theta^{*}_{1}\hspace{.3cm}....& \theta^{*}_{n} \\
 
 \end{bmatrix}
 \)$ = $ \(
 \begin{bmatrix}
 \alpha_{1}\hspace{.3cm}....& \alpha_{n}  \\
 \theta_{1}\hspace{.3cm}....& \theta^{*}_{n} \\
 
 \end{bmatrix}
 \)$T^{(i,j)}_{w}.$  Then, we have $ X_{n:n}\geq_{st}X^{*}_{n:n} .$
 \item Assume that \(
 \begin{bmatrix}
 \alpha_{1}\hspace{.3cm}....&\alpha_{n}  \\
 \theta_{1} \hspace{.3cm}....& \theta_{n} \\
 
 \end{bmatrix}
 \) $ \in \mathcal{S}_{n} ,$ and
 
 \(
 \begin{bmatrix}
 \alpha^{*}_{1}\hspace{.3cm}....& \alpha^{*}_{n}  \\
 \theta^{*}_{1}\hspace{.3cm}....& \theta^{*}_{n} \\
 
 \end{bmatrix}
 \)$ = $ \(
 \begin{bmatrix}
 \alpha_{1}\hspace{.3cm}....& \alpha_{n}  \\
 \theta_{1}\hspace{.3cm}....& \theta^{*}_{n} \\
 
 \end{bmatrix}
 \)$T^{(i,j)}_{w}.$  Then, we have $ X_{1:n}\leq_{st}X^{*}_{1:n} .$
 \end{enumerate}
 Where $ T^{(i,j)}_{w}=wI_{n}+(1-w)\Pi_{n}^{(i,j)}, $ $ \Pi_{n}^{(i,j)} $ is a permutation matrix of order $ n $ that just interchanges the coordinates $ i $ and $ j .$
 \begin{proof} 

 \begin{enumerate} 
 \item  Since $ T^{(i,j)}_{w}=wI_{n}+(1-w)\Pi_{n}^{(i,j)}. $  Here $ \Pi_{n}^{(i,j)} $ is a permutation matrix of order $ n $ that just interchanges the coordinates $ i $ and $ j ,$ so, in this case we have $ X_{k}$ and $ X^{*}_{k} $ to follow the same distribution, that is $ \alpha_{k}=\alpha^{*}_{k} $ and  $ \theta_{k}=\theta^{*}_{k},\hspace{.25cm} \forall  k\neq i,j. $ Then the theorem follows immediately from  Theorem \ref{th12}.
 \item The proof of (2) is similar to the proof of the part $ (1) $ and is therefore skipped.
 \end{enumerate}
 \end{proof}
 
 \label{th13}
 \end{theorem}

 It can be easily proved by the induction method that the product of a finite number of $ T $-transformation matrices with the same structure is also a $ T $-transformation matrix, and the structure is similar to the elements.  Following this fact, we obtain the corollary from Theorem \ref{th13}.
 
 \begin{corollary}
 Let $ X_{1},...., X_{n} $ and  $ X^{*}_{1},...., X^{*}_{n} $ be two sets of non-negative independent random variables such that $ X_{i}\sim EG_{2}(\theta_{i}, \phi, \alpha_{i}),$ and  $ X^{*}_{i}\sim EG_{2}(\theta_{i}^{*}, \phi, \alpha_{i}^{*})$ for $ i=1,...,n$.
 \begin{enumerate}
 \item  Assume that \(
 \begin{bmatrix}
 \alpha_{1}\hspace{.3cm}....&\alpha_{n}  \\
 \theta_{1} \hspace{.3cm}....& \theta_{n} \\
 
 \end{bmatrix}
 \) $ \in \mathcal{T}_{n} ,$ and
 
 \(
 \begin{bmatrix}
 \alpha^{*}_{1}\hspace{.3cm}....& \alpha^{*}_{n}  \\
 \theta^{*}_{1}\hspace{.3cm}....& \theta^{*}_{n} \\
 
 \end{bmatrix}
 \)$ = $ \(
 \begin{bmatrix}
 \alpha_{1}\hspace{.3cm}....& \alpha_{n}  \\
 \theta_{1}\hspace{.3cm}....& \theta_{n} \\
 
 \end{bmatrix}
 \)$T^{(i,j)}_{w_{1}}...T^{(i,j)}_{w_{k}}.$ \\
 
 Where $ T^{(i,j)}_{w_{i}},\hspace{.2cm}i=1,2,....,k$ have the same structure. Then, we have $ X_{n:n}\geq_{st}X^{*}_{n:n} .$
 \item Assume that \(
 \begin{bmatrix}
 \alpha_{1}\hspace{.3cm}....&\alpha_{n}  \\
 \theta_{1} \hspace{.3cm}....& \theta_{n} \\
 
 \end{bmatrix}
 \) $ \in \mathcal{S}_{n} ,$ and
 
 \(
 \begin{bmatrix}
 \alpha^{*}_{1}\hspace{.3cm}....& \alpha^{*}_{n}  \\
 \theta^{*}_{1}\hspace{.3cm}....& \theta^{*}_{n} \\
 
 \end{bmatrix}
 \)$ = $ \(
 \begin{bmatrix}
 \alpha_{1}\hspace{.3cm}....& \alpha_{n}  \\
 \theta_{1}\hspace{.3cm}....& \theta_{n} \\
 
 \end{bmatrix}
 \)$T^{(i,j)}_{w_{1}}...T^{(i,j)}_{w_{k}}.$ \\
 Where $ T^{(i,j)}_{w_{i}},\hspace{.2cm}i=1,2,...,k$ have the same structure. Then, we have $ X_{1:n}\leq_{st}X^{*}_{1:n} .$
 \end{enumerate}
 
 \label{corr14}
 \end{corollary}

 As we know, the finite product of $ T $-transformation matrices with different structures may not be a $ T $-transformation matrix. It is interesting to check whether the result in Corollary \ref{corr14} may still hold if the matrices $T_{w_{i}},\hspace{.1cm}i=1,2,...,k  $ have different structures. The following result gives an answer.
 
 \begin{theorem}
 Let $ X_{1},...., X_{n} $ and  $ X^{*}_{1},...., X^{*}_{n} $ be two sets of non-negative independent random variables such that $ X_{i}\sim EG_{2}(\theta_{i}, \phi, \alpha_{i}),$ and  $ X^{*}_{i}\sim EG_{2}(\theta_{i}^{*}, \phi, \alpha_{i}^{*})$ for $ i=1,...,n$.
 \begin{enumerate}
 \item  Assume that \(
 \begin{bmatrix}
 \alpha_{1}\hspace{.3cm}....&\alpha_{n}  \\
 \theta_{1} \hspace{.3cm}....& \theta_{n} \\
 
 \end{bmatrix}
 \) $ \in \mathcal{T}_{n} ,$ \par
 
 \(
 \begin{bmatrix}
 \alpha_{1}\hspace{.3cm}....& \alpha_{n}  \\
 \theta_{1}\hspace{.3cm}....& \theta_{n} \\
 
 \end{bmatrix}
 \)$T_{w_{1}}...T_{w_{i}}\in \mathcal{T}_{n}$\hspace{.3cm} for $ i=1,....,k-1,(k\geq2) .$  If\par
 
 \(
 \begin{bmatrix}
 \alpha^{*}_{1}\hspace{.3cm}....& \alpha^{*}_{n}  \\
 \theta^{*}_{1}\hspace{.3cm}....& \theta^{*}_{n} \\
 
 \end{bmatrix}
 \)$ = $ \(
 \begin{bmatrix}
 \alpha_{1}\hspace{.3cm}....& \alpha_{n}  \\
 \theta_{1}\hspace{.3cm}....& \theta_{n} \\
 
 \end{bmatrix}
 \)$T_{w_{1}}...T_{w_{k}},$ \\
 
 then, we have $ X_{n:n}\geq_{st}X^{*}_{n:n} .$
 \item Assume that \(
 \begin{bmatrix}
 \alpha_{1}\hspace{.3cm}....&\alpha_{n}  \\
 \theta_{1} \hspace{.3cm}....& \theta_{n} \\
 
 \end{bmatrix}
 \) $ \in \mathcal{S}_{n} ,$ \par
 
 \(
 \begin{bmatrix}
 \alpha_{1}\hspace{.3cm}....& \alpha_{n}  \\
 \theta_{1}\hspace{.3cm}....& \theta_{n} \\
 
 \end{bmatrix}
 \)$T_{w_{1}}...T_{w_{i}}\in S_{n}$\hspace{.3cm} for $ i=1,....,k-1, (k\geq2) .$ If\par
 
 \(
 \begin{bmatrix}
 \alpha^{*}_{1}\hspace{.3cm}....& \alpha^{*}_{n}  \\
 \theta^{*}_{1}\hspace{.3cm}....& \theta^{*}_{n} \\
 
 \end{bmatrix}
 \)$ = $ \(
 \begin{bmatrix}
 \alpha_{1}\hspace{.3cm}....& \alpha_{n}  \\
 \theta_{1}\hspace{.3cm}....& \theta_{n} \\
 
 \end{bmatrix}
 \)$T_{w_{1}}...T_{w_{k}},$ \\
 
 then, we have $ X_{1:n}\leq_{st}X^{*}_{1:n} .$
 \end{enumerate}
 \begin{proof}
 Set
 \(
 \begin{bmatrix}
 \alpha^{(j)}_{1}\hspace{.3cm}....& \alpha^{(j)}_{n}  \\
 \theta^{(j)}_{1}\hspace{.3cm}....& \theta^{(j)}_{n} \\
 
 \end{bmatrix}
 \)=
 \(
 \begin{bmatrix}
 \alpha_{1}\hspace{.3cm}....& \alpha_{n}  \\
 \theta_{1}\hspace{.3cm}....& \theta_{n} \\
 
 \end{bmatrix}
 \)$T_{w_{1}}...T_{w_{j}}\in \mathcal{T}_{n}$ for  $ j=1,....,k-1 .$\par
 Let $Y^{(j)}_{1},....,Y^{(j)}_{n},$\hspace{.1cm}  $ j=1,....,k-1 ,$ be sets of independent random variables such that $Y^{(j)}_{i}\sim EG_{2}(\theta_{i}^{(j)}, \phi, \alpha_{i}^{(j)}),$\hspace{.2cm} $ i=1,..,n, $ and $j=1,..,k-1.$
 \begin{enumerate}
 \item From the assumption of the theorem, it follows that 
 \begin{center}
 \(
 \begin{bmatrix}
 \alpha^{(j)}_{1}\hspace{.3cm}....& \alpha^{(j)}_{n}  \\
 \theta^{(j)}_{1}\hspace{.3cm}....& \theta^{(j)}_{n} \\
 
 \end{bmatrix}
 \)$ \in \mathcal{T}_{n} $
 \end{center} for $j=1,..,k-1, (k\geq 2).$ Using the result of Theorem \ref{th13} and these observations, it follows that  $ X_{n:n}\geq_{st}Y^{(1)}_{n:n}\geq_{st}....\geq_{st}Y^{(k-1)}_{n:n}\geq_{st} X^{*}_{n:n},$ which completes the proof of part (1).
 \item The proof is similar to the proof of part (1) and is therefore skipped here for the sake of brevity.
 \end{enumerate}
 \end{proof}
 
 \label{th15}
 \end{theorem}

\section{Discussion and Conclusion}
\label{con}

In this paper, we have presented a comprehensive stochastic ordering analysis of the Exponentiated Gumbel Type-II ($EG_2$) distribution, focusing on how parameter heterogeneity influences the reliability of series and parallel systems. By employing the theory of majorization and multivariate chain majorization, we established several ordering results that provide both theoretical and practical insights.

Our findings demonstrate that the $EG_2$ model exhibits a rich and complex structural behavior that distinguishes it from simpler exponentiated families. First, for series systems (lowest order statistics), we derive the usual stochastic and likelihood ratio orders, which further imply failure rate, mean residual life, and increasing convex orderings.  We also establish that the model admits a proportional failure rate structure. In addition, we establish star and dispersive orderings, which provide further insight into variability and aging-related characteristics of system lifetimes.
Second, for parallel systems (highest order statistics), the results reveal a more nuanced behavior. While we establish the reversed failure rate order, we also identify a \emph{non-comparability} phenomenon. Specifically, unlike many classical models, the highest order statistics of $EG_2$ variables are not necessarily comparable under the usual stochastic order when the inner shape parameters are heterogeneous. This highlights the sensitivity of parallel systems to component-level characteristics and cautions against extending conclusions from simpler lifetime models.
Finally, by extending our analysis to multivariate chain majorization, we provide a framework for comparing systems in which scale and outer shape parameters vary simultaneously, allowing for more general reliability comparisons.

The stochastic ordering results in this work are derived under structural assumptions on the model parameters, such as independence of components, equality of certain parameters across components and systems, and membership of parameter vectors and matrices in the sets $ \mathcal D_n$, $ \mathcal S_n$, and $ \mathcal T_n$. These conditions are imposed to ensure the Schur-convexity or Schur-concavity of the relevant reliability functions and to make the majorization arguments tractable. Therefore, these assumptions should be viewed as sufficient but not necessary for the stated orderings.
Extending the present results to more general settings is an interesting direction for future research. Another direction for future work is to investigate stochastic comparisons under dependent component lifetimes, for example, using copula-based models, as well as more complex system configurations such as $k$-out-of-$n$ systems.

\section*{Acknowledgment}
The author sincerely thanks the reviewers for their valuable comments and suggestions, which have significantly improved this manuscript. The author also expresses gratitude to Dr. Nitin Gupta from the Department of Mathematics, IIT Kharagpur, India, for helpful discussions and support during the preparation of this work.

\section{Declarations} 
    \subsection*{Funding} This research was not funded by any specific source.

    \subsection*{Data availability} Not Applicable.
\subsection*{Conflict of interest}The author declares no conflict of interests/competing interests.

\end{document}